\documentclass[12pt]{article}
\usepackage{pudlatex,fullpage}
\usepackage{amsmath}
\usepackage[all,cmtip]{xy}

\hi

\title{The canonical pairs of bounded depth Frege systems}

\author{Pavel Pudl\'ak
\thanks{The author is supported by the project EPAC, funded by the Grant Agency of the Czech Republic under the grant agreement no. 19-27871X,  and the institute grant RVO: 67985840. Part of this work was done when the author was supported by the ERC Advanced Grant  339691 (FEALORA).}}

\maketitle

\begin{abstract}
The canonical pair of a proof system $P$ is the pair of disjoint NP sets where one set is the set of all satisfiable CNF formulas and the other is the set of CNF formulas that have $P$-proofs bounded by some polynomial. We give a combinatorial characterization of the canonical pairs of depth~$d$ Frege systems. Our characterization is based on certain games, introduced in this article, that are parametrized by a number~$k$, also called the depth. We show that the canonical pair of a depth~$d$ Frege system is polynomially equivalent to the pair $(A_{d+2},B_{d+2})$ where $A_{d+2}$ (respectively, $B_{d+1}$) are depth {$d+1$} games in which Player~I (Player II)  has a positional winning strategy. Although this characterization is stated in terms of games, we will show that these combinatorial structures can be viewed as generalizations of monotone Boolean circuits. In particular, depth~1 games are essentially monotone Boolean circuits. Thus we get a generalization of the monotone feasible interpolation for Resolution, which is a property that enables one to reduce the task of proving lower bounds on the size of refutations to lower bounds on the size of monotone Boolean circuits. However, we do not have a method yet for proving lower bounds on the size of depth~$d$ games for $d>1$.
\end{abstract}

\section{Introduction}

There are two basic problems associated with every propositional proof system $P$:
\ben
\item Given a CNF formula $\phi$, decide whether $\phi$ is satisfiable or has a short $P$ refutation, provided that we know that one of these statements is true.
\item Let two CNF formulas $\phi$ and $\psi$ with disjoint sets of variables and a refutation $\pi$ of $\phi\wedge\psi$ be given and assume that one of the formulas is satisfiable. Decide which.
\een
These problems are formalized by defining pairs of disjoint NP sets. The pair of the first problem is called the \emph{canonical pair} of $P$ (introduced by Razborov in~\cite{razb-NP}), the pair of the second one is called the \emph{interpolation pair} (introduced in~\cite{pudlak:03}). We will see that these two problems are tightly connected and thus one can focus only on one of them. Since interpolation pairs are not so well-known as canonical pairs, we use this name in the title, but in fact, we will study interpolation pairs.

We conjecture that the hardness of these problems increases with the strength of the proof systems, where we compare hardness of disjoint NP pairs by polynomial reductions. We do not have means to prove that one pair is stronger than another, because if P=NP all are interreducible. So we have only two possibilities what to do. First, we can try and find mathematical principles equivalent to the facts that these pairs are disjoint. If the principles seem of increasing strength we can view it as evidence of increasing hardness of these pairs. Second, we can consider monotone versions of these problems. Since we do have lower bounds on monotone Boolean circuits and some other monotone computational models (in contrast to the desperate state of the affairs with general Boolean circuits), there is some chance that we can prove separation at least with respect to monotone reductions.

It is natural to start with the weakest systems. Prior to this work, combinatorial characterization of the canonical pair was known only for the Resolution system~\cite{BPT}. Bounded depth Frege is a well-studied hierarchy of proof systems above Resolution. In this article we will present combinatorial characterizations of interpolation pairs for all levels of this hierarchy, which also gives characterizations of canonical pairs. 

Our characterization is based on certain two player games which we will briefly explain now and give precise definition later. Two players alternate in writing symbols on a finite tape. They start on one end, say the left one, and proceed to the other. When they reach the end they either stop, if the game has only one round, or they reverse the direction and go back. They may reverse direction $(k-1)$-times if the depth parameter is $k$. What is a legal move only depends on the symbol in the current square and the next one. We define \emph{positional} strategies and show that given a positional strategy one can decide in polynomial time whether it is a winning strategy for the particular player. This enables us to define an NP pair for every depth $k\geq 1$ and characterize the interpolation pair of depth-$d$ Frege systems by games of depth $d+1$. The canonical pair of depth $d$ Frege systems is polynomially equivalent to the interpolation pair of depth $d+1$ Frege systems.

One can view a game of depth $k+1$ as follows. In the first round the players alternate to define a string of symbols that determine a game of depth $k$ that is played after the first round. This suggests the intuition that it should be harder to decide who has a (positional) winning strategy in a game of depth $k+1$: we cannot use an oracle for games of depth $k$ because the game of depth $k$ is yet to be determined by playing the first round.

Our result can also be viewed as a contribution to the line of research that studies monotone computation models. We will show that one can interpret our games, more precisely game schemas, as monotone computation models generalizing monotone Boolean circuits. It is not clear how difficult it may be to prove lower bounds on these models, but if we could do it, we may be able to solve an important problem about bounded depth Frege systems. Impagliazzo and Kraj\'{\i}\v{c}ek~\cite{impagliazzo-krajicek} proved that depth $d$ Frege systems cannot polynomially simulate depth $d+1$ Frege systems w.r.t. refuting CNFs. However their lower bound is only mildly superpolynomial, while we believe there should be exponential separation. The tautologies, or rather contradictions, based on games are candidates for exponential separations. 
%% Unlike other candidate contradiction, they can be stated as 3-CNFs (???).

A reader familiar with results in proof complexity and Bounded Arithmetic will recognize several connections between this work and previous ones. The Symmetric Calculus is inspired by the calculus invented by Skelley and Thapen~\cite{ST}. The \emph{Game Induction} principles introduced in~\cite{ST} are also very much related to our bounded depth games. It is possible that some arguments from that article could be used for proving our result, similarly as one might use Skelley-Thapen's calculus instead of our Symmetric Calculus. A game similar to ours appeared in an article of Ko{\l}odziejczyk, Nguyen and Thapen~\cite{KNT} (Lemma~10). We will also see in Section~\ref{s8} that the \emph{point-line game}, introduced an article of Beckmann, Pudl\'ak, and Thapen~\cite{BPT}, can be viewed as a version of our depth-2 game. There are certainly more such connections than those mentioned above.

%\bigskip
This article is essentially a proof of a single theorem, Theorem~\ref{t-main}, plus some observations. We start by recalling the definitions of bounded depth sequent calculi (that are used instead of bounded depth Frege calculi), canonical and interpolation pairs and stating some basic facts. Then we define the \emph{Symmetric Calculus}. This calculus, more precisely its bounded depth version, has been designed for proving our theorem, but it may be of independent interest. The basic idea is due to Skelley and Thapen~\cite{ST}, but our calculus differs in several particulars. We show that the bounded version of Symmetric Calculus is polynomially equivalent to the standard formalization by the sequent calculus. In Section~\ref{s4} we define the games used in the characterization. In Sections~\ref{s5} and~\ref{s6} we construct the reductions. In Section~\ref{s7} we prove a stronger version of our main theorem. In Section~\ref{s8} we will have a closer look at games of depth~1 and~2. We conclude the article with some open problems. At the end there is a short appendix in which we mention a connection to Bounded Arithmetic, which is important, but not used in this article, and explain one technical point from the simulation of bounded depth calculi by the bounded depth Symmetric Calculi.

\paragraph{Acknowledgment.} I am grateful to Emil Je\v{r}\'abek and Jan Kraj\'{\i}\v{c}ek for their comments on the draft of this article and especially to Neil Thapen for reading the whole manuscript an pointing to incomplete, or unclear parts.

\section{Basic notions}\label{s2}

In this section we recall some concepts and results from proof complexity that we will use.

\subsection{Bounded depth sequent calculus}

Classical propositional logic can be formalized by various types of calculi, which may have different power. We compare calculi by how efficiently they can prove tautologies. A tautology that can only be proved by exponentially long proofs in one system may have polynomial size proofs in another. The standard formalization of propositional logic is based on axioms and derivation rules. These calculi are called \emph{Frege calculi}. They are equivalent, from the point of view of efficiency, to the \emph{sequent calculus}. 

In this article we are interested in restricted versions of these calculi where the depth of formulas is bounded by a constant. To this end we will restrict our language to the De Morgan basis $\neg,\vee,\wedge$. The \emph{depth} of a formula is defined to be the number of alternations of $\neg,\vee,\wedge$, where negations at variables are not counted. Thus variables and negated variables, called \emph{literals}, have depth~0, conjunctions of variables and negated variables have depth~1, disjunctions of variables and negated variables, called\emph{clauses}, have depth~1, CNFs and DNFs have depth~2, etc.

For formalizing reasoning with bounded depth formulas, the sequent calculus is more convenient than Frege calculi. We define \emph{depth~$d$ sequent calculus} to be the standard sequent calculus restricted to the De Morgan basis and formulas of depth at most $d$. In the sequel we will only use bounded depth sequent calculi, but the reader should keep in mind that they are equivalent to bounded depth Frege calculi.

Another useful convention is to use \emph{refutations} instead of derivations. Given a DNF tautology $\tau$, we take the CNF contradictory formula $\sigma$ obtained by the dualization of $\tau$ and prove contradiction from $\sigma$ (contradiction is represented by the empty sequent). Since a CNF formula can be represented by a set of clauses and clauses can be represented by sequents consisting of literals only, we can use even the depth~0 calculus. The depth~0 calculus is essentially the \emph{Resolution} system, because the only non-structural rule that can be used is cut with a literal as the cut formula.

We will shortly introduce yet another calculus for reasoning with bounded depth formulas, the \emph{Symmetric Calculus}. In this calculus, $\Pi_{k+2}$ proofs correspond to proofs in the depth~$k$ sequent calculus.

\subsection{Polynomial simulations and disjoint NP pairs of propositional proof systems}

We say that a proof system \emph{$P$ polynomially simulates a proof system $Q$} if there exists a polynomial time algorithm that from a given $P$-proof (or refutation) $\pi$ of a formula $\phi$, constructs a $Q$-proof (or refutation) of $\phi$. We say that \emph{$P$ and $Q$ are polynomially equivalent} if they polynomially simulate each other.

Let $(A,B)$ and $(C,D)$ be two pairs of disjoint NP sets. We say that $(A,B)$ is \emph{polynomially reducible} to $(C,D)$ if  there exists a polynomial time algorithm that maps $A$ to $C$ and $B$ to~$D$.

The \emph{canonical pair of a proof system $P$} (defined in~\cite{razb-NP}) is the pair of disjoint NP-sets $(A,B)$ where
\[\ba{l}
A:=\{(\phi,0^m)\ |\ \phi\mbox{ satisfiable }\},\\
B:=\{(\phi,0^m)\ |\ \phi\mbox{ has a $P$ refutation of size }m\}.
\ea\]
The string of zeros $0^m$ of length $m$ is a padding that enables us to consider proofs of arbitrary length. In all natural proof systems we can replace this padding by padding the formulas with trivially satisfiable clauses. Then we can define $A$ to be satisfiable formulas and $B$ to be formulas $\phi$ that have $P$-refutations of length $|\phi|^2$.

The \emph{interpolation pair of a proof system $P$} (defined in~\cite{pudlak:03}) is the pair of disjoint NP-sets $(A,B)$ where
\[\ba{l}
A:=\{(\phi,\psi,\pi)\in\Delta\ |\ \phi\mbox{ satisfiable }\},\\
B:=\{(\phi,\psi,\pi)\in\Delta\ |\ \psi\mbox{ satisfiable }\},
\ea\]
where $\Delta$ is the set of triples $(\phi,\psi,\pi)$ such that $\phi$ and $\psi$ are formulas with disjoint sets of variables and $\pi$ is a $P$-refutation of $\phi\wedge\psi$. 

In these definitions we have not specified the class of Boolean formulas. Propositional proof systems may use restricted classes of formulas and if they use different classes, then polynomial simulation does not make sense. A natural minimal requirement is that a proof system is complete with respect to refutations of unstatisfiable CNF formulas. Therefore \emph{we restrict the above definitions to CNF formulas.} One can show, under very mild assumptions about the proof systems, that such restricted pairs are polynomially equivalent to the pairs defined above. Hence we do not lose information about the complexity of these pairs if we focus on CNFs.

The basic facts about these concepts are:
\ben
\item {\it the interpolation pair of $P$ is polynomially reducible to the canonical pair of $P$;
\item if $P$ polynomially simulates $Q$, then the canonical (respectively, interpolation) pair of $Q$ is polynomially reducible to the canonical (interpolation) pair of $P$.}
\een

For bounded depth sequent calculi, we have the following important fact.

\bpr[\cite{BPT}, Proposition 1.4]\label{p2.1}
For $k\geq 0$, the canonical pair of depth~k sequent calculus is polynomially equivalent to the interpolation pair of depth~(k+1) sequent calculus. 
\epr
Therefore it suffices to characterize the interpolation pairs.

With each disjoint NP pair $(A,B)$, there is an associated \emph{separation  problem:} given the promise that $x\in A\cup B$, how difficult is to decide whether $x\in A$ or $x\in B$? For the interpolation pair of depth~1 sequent calculus, which is the same as Resolution, the problem is decidable in polynomial time. For higher systems we do not know if the separation is solvable in polynomial time. By improving some previous results, Bonet et al.~\cite{bonet} proved that assuming factoring Blum integers or computing the Diffie-Helman function is sufficiently hard, the separation problem for the interpolation pairs is not polynomially solvable for all depth~$d$ sequent calculi starting from some small $d_0$. Less convincing evidence of the hardness are results showing that the decision of who has a winning strategy in certain combinatorial games can be reduced to the interpolation pairs of some small depth~$d$ sequent calculi~\cite{huang-pitassi,atserias-maneva,BPT}. In particular, the decision problem for \emph{parity games} can be reduced to the canonical pair of depth~0 sequent calculus, i.e., Resolution, and the decision problem for \emph{simple stochastic games} can be reduced to the canonical pair of depth~1 sequent calculus. The decision problem for parity games is solvable in quasipolynomial time, but the author of this article thinks that the canonical pair of Resolution is harder. This belief is supported by the recent result of Atserias and M\"uller~\cite{atserias-muller} that the proof search in the Resolution system is NP-hard.

\subsection{Feasible interpolation}
%\subsection{Monotone interpolation and monotone reductions}

\paragraph{The feasible interpolation property.} We say that a proof system $P$ has the \emph{feasible interpolation property} if there exists a polynomial time algorithm $A$ such that given a $P$-proof $D$ of a formula of the form
\[
\phi(\bar z,\bar x)\vee\psi(\bar z,\bar y),
\]
where all common variables of the formulas are in the string $\bar z$, and given an assignment $\bar a\in\{0,1\}^n$ to variables $\bar z$, the following holds true:
\bel{e-int}
\baa{l}
\mbox{ if }\phi(\bar a,\bar x) \mbox{ is satisfiable, then }A(D,\bar a)=0,\\
\mbox{ if }\psi(\bar a,\bar y) \mbox{ is satisfiable, then }A(D,\bar a)=1.
\ea
\ee
It follows that for every $D$, there exists a \emph{Boolean circuit} $C(\bar z)$ whose size is polynomial in the size of $D$ and such that condition~(\ref{e-int}) is satisfied with $A(D,\bar a)$ replaced with $C(\bar a)$.
For natural proof systems,%
\footnote{We only need that given a $P$ proof of $\alpha(\bar z,\bar x)$, one can construct in polynomial time a proof of $\alpha(\bar a,\bar x)$ for every assignment $\bar a$ to $\bar z$.}
$P$ has the feasible interpolation property iff the interpolation pair of $P$ is separable by a polynomial time algorithm. 

\paragraph{The monotone feasible interpolation property.} Often one can show that there exists a \emph{monotone} Boolean circuit with the above properties. For this, it is necessary to ensure that the two sets can be separated by monotone functions, which is done by assuming that all common variables $\bar z$ occur only negatively in $\phi$, or all occur positively in $\phi$ (or both). If for this kind of formulas, there exist polynomial size monotone circuits with property~(\ref{e-int}), then we say that the proof system has the  \emph{monotone feasible property.} 

It is well-known that the Resolution proof system has both the feasible interpolation property and the monotone feasible property. There are a few more natural proof systems that have the feasible interpolation property, some also have the monotone feasible property, see~\cite{krajicek}, Chapter~17. Our Theorem~\ref{t7.1} can be viewed as a generalization of the result for Resolution to stronger fragments of the propositional sequent calculus.

%% \paragraph{Monotone reductions.} If in the formula $\phi(\bar z,\bar x)\vee\psi(\bar z,\bar y)$ common variables $\bar z$ occur negatively in $\phi$ and  positively in $\phi$, then the formula defines in natural way a partial monotone Boolean function $f$. Monotone interpolation means that this function can be defined by a monotone Boolean circuits whose size is polynomially bounded by the size of the proof. In this article we will present characterizations that can be viewed as kinds of monotone computational devices. To 

%% \bigskip
%% TBD

\section{The Symmetric Calculus}\label{s3}

This calculus is specially designed for the proof of the main theorem. The starting point was the calculus of Skelley and Thapen~\cite{ST}. In their calculus the eigenformulas of non-structural rules are just literals. Thus instead of the general cut rule they use the resolution rule, i.e., cut with a literal as the eigenformula, and the rule for conjunction introduction only allows a  conjunction to be extended by a literal. The fact that such a calculus can polynomially simulate bounded depth sequent calculus is a remarkable discovery, because it is well-known that if one restricts cuts in the sequent calculus to depth $d$ formulas, then the resulting system has only the power of depth $d$ sequent calculus (for refuting CNFs). The reason why restriction to literals does not limit the Skelley-Thapen calculus is that the calculus uses deep inferences.

In the Symmetric Calculus, unlike the Skelley-Thapen calculus, all rules can be applied as deep inferences. In order to achieve symmetry of the rules, we have replaced conjunction introduction by dual resolution. Furthermore, the proofs do not have the traditional structure where the set of initial formulas is gradually extended by derived formulas. In the Symmetric Calculus a proof has a linear structure---a sequence of formulas such that the next formula follows only from the previous formula.

In the bounded version there are further structural restrictions. Let us repeat that this is only because we need to make a connection with certain games.

\paragraph{Definition of the Symmetric Calculus.}
The language of the Symmetric Calculus consists of $\vee,\wedge,\bot,\top$,  and literals $x_i,\neg{x}_i$. Negations are allowed only at literals. Given a literal $p$, we denote by $\neg p$ its dual.

The calculus is based on \emph{deep inferences}, which means that one can replace a \emph{subformula} by a formula allowed by a rule. An application of a rule $\frac BC$ is a substitution
\[
\frac{A[\dots B\dots]}{A[\dots C\dots]}
\]
In the Symmetric Calculus every rule has one assumption and one conclusion, so the calculus is a term rewriting system.
% We only allow negations at literals, hence we have $B\models C$ ($C$ follows semantically from $B$) for every rule $\frac BC$ of the calculus.

A proof of $\Phi\vdash\Psi$ is a sequence of formulas $\Phi=\Phi_1\dts\Phi_m=\Psi$ where $\Phi_{i+1}$ follows from $\Phi_i$ by an application of a deduction rule.

\paragraph{The rules of the calculus:}~

\medskip
{\it commutativity and associativity of $\vee$ and $\wedge$,}%
\footnote{Emil Je\v{r}\'abek has observed that associativity is redundant and if we add also the other two versions of weakenings, $\frac{B}{A\vee B}$, $\frac{A\wedge B}{B}$, then also commutativity will be redundant.}

\medskip
{\it contraction/cloning}
\[
\frac{A\vee A}{A}\hskip 4cm \frac{A}{A\wedge A}
\]

\medskip
{\it $\bot$-elimination / $\top$-introduction}
\[
\frac{A\vee \bot}{A}\hskip 4cm \frac{A}{A\wedge \top}
\]

\medskip
{\it weakenings}
\[
\frac{A}{A\vee B}\hskip 4cm \frac{A\wedge B}{A}
\]

\medskip
{\it dual resolution/resolution}
\[
\frac{A\wedge\top\wedge B}{(A\wedge p)\vee(B\wedge\neg p)}\hskip 2cm
\frac{(A\vee p)\wedge(B\vee\neg p)}{A\vee\bot\vee B}
\]
In the last two rules $A$ or $B$ or both formulas may be not present; e.g., {\Large $\frac\top{~p\vee\neg p}$} is considered to be an instance of dual resolution.

\bigskip
The reason for calling this calculus ``symmetric'' is that each rule has its dual. Hence, given a proof of $\Phi\vdash\Psi$ we can obtain a proof of $\neg\Psi\vdash\neg\Phi$, by inverting the order of formulas and replacing connectives, truth constants, and literals by their duals. Here and also in the sequel, $\neg\Phi$ denotes the dual of the formula $\Phi$. More importantly, the symmetry allows us to cut the case analysis to one half.

Several comments about the rules are in order. 
The reason for using the rules for the truth constants in this form instead of the standard ones
\[
\frac{\bot}{A}\hskip 2cm \frac{A}{\top}
\]
is purely technical. The weakenings are also called ``disjunction introduction'' and ``conjunction elimination''. We prefer to view them as weakenings, because we do not have rules for disjunction elimination and conjunction introduction.%
\footnote{Another reason is that the rule of weakening can be omitted in sequent calculi if we use two versions of every other rule: one in which the premises are ``consumed'' and one in which they stay. This might be possible also in this calculus, but we have not investigated this possibility.}
 The truth constants in resolution and dual resolution rules can, clearly, be omitted when at least one of the formulas $A$ or $B$ is present, but, again, for technical reasons, we prefer to keep the constants also when the context is nonempty.

\subsection{Cuts and dual cuts}\label{s-4.1}

Our aim now is to prove that dual cuts and cuts can be simulated. These derived rules are generalizations of the rules of dual resolution and resolution where one can use arbitrary formulas $C$ instead of literals:
\[
\frac{A\wedge\top\wedge B}{(A\wedge C)\vee(B\wedge\neg C)}\hskip 2cm
\frac{(A\vee C)\wedge(B\vee\neg C)}{A\vee\bot\vee B}
\]
where $\top$ and $\bot$ can be omitted if at least one of the formulas $A$, $B$ is present. Recall that $\neg C$ denotes the formula obtained from $C$ by replacing the connectives and literals by their duals. 
% We will only show that the simulation is polynomial for constant logical depth. 

We will start with a simple fact.

\bfa
The distributivity of $\vee$ over $\wedge$ can be polynomially simulated.
\efa
\bprf
\[

\ba{ccl}
&A\vee(B\wedge C)&\\
\cline{1-2}
&(A\vee(B\wedge C))\wedge(A\vee(B\wedge C))&\mbox{\quad by cloning}\\
\cline{1-2}
&(A\vee B)\wedge(A\vee C)&\mbox{\quad by weakenings}
\ea\]
\eprf
The proof above explains what we mean by \emph{polynomial simulation:} if $\psi$ is a formula obtained from $\phi$ by replacing an occurrence of $\alpha\vee(\beta\wedge\gamma)$ with $(\alpha\vee\beta)\wedge(\alpha\vee\gamma)$, then one can construct in polynomial time a derivation of $\psi$ from $\phi$ in the Symmetric Calculus.

Using the symmetry of our calculus we immediately get that
\[
\frac{(A\wedge B)\vee(A\wedge C)}{A\wedge(B\vee C)}
\]
also can be polynomially simulated.
The distributivity of $\wedge$ over $\vee$ can be simulated too, but  to prove it, we first need to simulate cuts and dual cuts.

\bll{l-1}
Dual cuts and cuts can be polynomially simulated.
\el
\bprf
We will describe a procedure that constructs a proof that simulates cut. Suppose we want to simulate the following deduction
\[
\frac{(A\vee C)\wedge(B\vee \neg C)}{A\vee B}.
\]
If $C$ is a literal, then this is just an application of the resolution rule. Now suppose that $C$ is $C_1\wedge C_2$ (the case of $C$ being $C_1\vee C_2$ will follow by symmetry). First we use distributivity to obtain
\[
(A\vee C_1)\wedge(A\vee C_2)\wedge(B\vee \neg C_1\vee \neg C_2).
\]
Then we apply the procedure recursively to the subformula $(A\vee C_1)\wedge(B\vee \neg C_1\vee \neg C_2)$ (after permuting $A\vee C_1$ and $A\vee C_2$). Thus we obtain
\[
(A\vee C_2)\wedge(B\vee \neg C_2).
\]
In this way we have reduced the problem to a smaller cut-formula $C_2$ and we can recursively call the procedure again.

To see that this gives a polynomial simulation it suffices to observe two facts:
\ben
\item the number of times the procedure calls itself is equal to the number of subformulas of~$C$;
\item each call of the procedure adds a term to the initial conjunction $(A\vee C)\wedge(B\vee \neg C)$ whose size is bounded by either the size of $A\vee C$ or $B\vee\neg C$. In the derivation above the added formula was first a clone of the formula $A\vee C$ and then it was weakened to $A\vee C_2$, while we consider $A\vee C_1$ to be only a weakened original $A\vee C$.
\een

%% \[
%% \def\arraystretch{1.5}
%% \begin{array}{lcl}
%% (1)&(A\vee C)\wedge(B\vee \neg C)&\\
%% \cline{2-2}
%% (2)&(A\vee \bigwedge_iC_i)\wedge(B\vee \bigvee_i\neg C_i)&\mbox{\quad the same formula}\\
%% \cline{2-2}
%% (3)&\bigwedge_i(A\vee C_i)\wedge(B\vee \bigvee_i\neg C_i)&\mbox{\quad by distributivity}\\
%% \cline{2-2}
%% (4)&A\vee B&\mbox{\quad by cuts with $C_i$s and contractions}
%% \end{array}\]

The case of dual cuts follows by symmetry. 
See also Appendix for a remark and an example.\eprf

% We do not know if this simulation is polynomial in general, i.e., without the restriction to $\Pi_k$ formula for a fixed $k$.

We can now show the simulation of the distributivity of $\wedge$ over $\vee$.

\bll{l-4.2}
The distributivity of $\wedge$ over $\vee$ with can be polynomially simulated.
\el
\bprf
\[
\baa{ccl}
&A\wedge(B\vee C)&\\
\cline{1-2}
&((A\wedge B)\vee \neg B)\wedge ((A\wedge C)\vee \neg C)\wedge(B\vee C))&\mbox{\quad by cloning $A$ and dual cuts}\\
\cline{1-2}
&(A\wedge B)\vee (A\wedge C)&\mbox{\quad by cuts with $B$ and $C$}
\ea\]
\eprf

By symmetry we get immediately
\bel{e-distr}
\frac{(A\vee B)\wedge(A\vee C)}{A\vee(B\wedge C)}
\ee
This is important because it enables us to simulate \emph{conjunction introduction} in the Sequent Calculus. Since we also have dual resolution (we can derive $p\vee\neg p$ by $\top$-introduction and dual cut), we have simulations of all rules of the Sequent Calculus as was formalized by Tait, cf.~\cite{buss-intro}. Thus we have shown:

\bpr
The Symmetric Calculus polynomially simulates the Sequent Calculus, hence also Frege calculi.
\epr

We will also need 
\bel{e-ext-con}
\frac{A\to(B\vee C)}{(A\wedge D)\to((B\wedge D)\vee C)}
\ee
If we write $\to$ in terms of $\vee,\neg$, this becomes
\[
\frac{\neg A\vee B\vee C}{\neg A\vee\neg D\vee(B\wedge D)\vee C}
\]
which is simply an application of dual cut.

We will show that these simulations also hold for bounded depth versions of the calculi; this will be more complicated.

%%%%%%%%%%%%%%%%%%%%%%%%%%%%%%%%%%%%%%%%%%%%%%%%%%%%%%%%%%%%%%%%%%

\subsection{The Bounded Depth Symmetric Calculus}\label{s4.2}

By \emph{logical depth} we mean the number of alternations of conjunctions and disjunctions. More precisely, literals are defined to be $\Sigma_0=\Pi_0$, conjunctions of literals are $\Pi_1$, disjunctions of literals are $\Sigma_1$ etc.

In the bounded depth calculus conjunctions and disjuctions are  operations with an arbitrary finite number of arguments $n=1,2,3,\dots$. The rules of commutativity and associativity are replaced by a general rule of permutation.%
\footnote{Why don't we use multisets? The reason is that we want to have correspondence between subformulas of a formula and subformulas in its successor in a proof, which is needed for defining games from proofs.}
The other rules of the Symmetric Calculus are applied to consecutive formulas in a possibly longer disjunction, or conjunction, and the derived formula is inserted on the position of the premise. For instance, an application of the contraction rule transforms a disjunction into a disjunction with one term less as follows:
\[
\frac{B_1\vee\dots\vee B_k\vee A\vee A\vee C_1\vee\dots\vee C_l}{B_1\vee\dots\vee B_k\vee A\vee C_1\vee\dots\vee C_l}
\]

\paragraph{Stratified formulas.} Furthermore, we require that $\wedge$ and $\vee$ alternate regularly in formulas as we go from the top connective to the bottom of the formula. Balanced formulas with $\wedge$ and $\vee$ alternating regularly will be called \emph{stratified}. More precisely, \emph{stratified $\Pi_k$ formulas} (\emph{$\Sigma_k$ formulas}) have the structure of rooted trees in which
\ben
\item every branch has length $k$ (measured by the number of edges),
\item on every branch $\wedge$ and $\vee$ alternate regularly starting with $\wedge$ (respectively with $\vee$),
\item the leaves are labeled by literals.
\een
We will use $\Pi^s_k$ and $\Sigma^s_k$ for the classes of stratified $\Pi_k$ and $\Sigma_k$ formulas. Note that formulas in classes $\Pi^s_k$ and $\Sigma^s_k$ have depth \emph{exactly}~$k$.

In order to represent formulas that are not in this form, we will use unary operations of $\wedge$ and $\vee$. We will use prefix notation for these unary conjunctions and disjunctions, e.g., $\wedge(A)$, or just $\wedge A$, while keeping infix notation for strings of formulas with at least two terms.

\ex{ 
The formula $(p\vee q)\wedge r$ can be represented by the stratified formula $(p\vee q)\wedge(\vee r)$. Note that in general there may be several different representations (see zipping and unzipping below).}

We will say that a formula $\phi$ is a \emph{legal subformula} of $\psi$ if it is a formula determined by a node $N$ in the tree representing $\psi$ in the following sense: the tree of $\phi$ is the entire tree below~$N$. 

%%%%% Note that the specification ``stratified'' means more than just that $\phi$ is stratified.

\ex{In the formula $(p\wedge q\wedge r)\vee(s\wedge t)$, the formula $p\wedge q\wedge r$ is a legal subformula, whereas $p\wedge q$ isn't.}

We will study $\Pi^s_k$ proofs of contradiction. In such a proof $A_1,A_2\dts A_n$, \emph{every formula $A_i$ must be $\Pi^s_k$}. E.g., if we are refuting a CNF formula $A$, we must pad it to a $\Pi^s_k$ formula and the final $\bot$ must be padded to level $\Pi^s_k$ too.

The deep inferences of the Symmetric Calculus may be applied only to legal subformulas. To this end we have to modify the rules. 
We will only define the modification of the rules in the left column of the list of rules in the previous subsection; the right column is done symmetrically.

\ben
\item Contraction means that we can replace two consecutive terms in a disjunction when they are equal.
\item Elimination of $\bot$ means that we remove it from the disjunction. This is not allowed if it is the only term in the disjunction (unlike in Resolution, in the Symmetric Calculus empty disjunctions are not used).
\item Weakening of a disjunction means inserting an arbitrary formula of appropriate logical complexity into the disjunction on arbitrary place.
\item Dual resolution means that we split a conjunction into the part before $\top$ and the part after, omit $\top$, and add a literal to the first part and the dual literal to the second:
\[
\frac{A_1\wedge\dots\wedge A_i\wedge\top\wedge A_{i+1}\wedge\dots\wedge A_n}{(A_1\wedge\dots\wedge A_i\wedge p)\vee( A_{i+1}\wedge\dots\wedge A_n\wedge\neg p)}.
\]
\een

\ex{Let $A,B,C\in\Sigma^s_k$. The following is \emph{not} legal application of the dual resolution in the bounded depth Symmetric Calculus
\[
\frac{A\wedge\top\wedge B\wedge C}{((A\wedge p)\vee(B\wedge\neg p))\wedge (\vee\wedge C)},
\]
even though the conclusion is stratified, because the rule is not applied to the entire conjunction. (For the sake of readability, we have omitted the padding of literals and~$\top$.)
}

We require that the rules be applied so that stratification is preserved. If contraction is applied to a disjunction in which there are only the two terms $A\vee A$, then the result is a unary disjunction $\vee(A)$; similarly for $\bot$-elimination. Weakening {\large $\frac{A}{A\vee B}$} can only be applied if $A$ is a part of disjunction, which may be just $\vee(A)$. Similarly for dual resolution, $A_1\wedge\dots\wedge A_i\wedge\top\wedge A_{i+1}\wedge\dots\wedge A_n$ must be a term in a disjunction, possibly the unique term in the disjunction. Thus we have
\bi
\item all rules in the left column can only be applied to disjunctions and
\item all rules in the right column can only be applied to conjunctions.
\ei

%% Consider, e.g., a stratified $\Pi^s_k$ formula 
%% $$(A\vee A)\wedge B$
%% where $A$ is $\Pi^s_{k-2}$.
%% After we contract $A\vee A$ to a $A$,  $A$ should be padded to a $\Sigma^s_{k-1}$ formula (by adding a unary $\vee$ on the top). This is ensured by the proper interpretation of how the rules are applied: \emph{we contract disjunction with two elements to a disjunction with one element}. 

%In the tree representation, contraction means that we delete one of the identical subtrees of a vertex labeled by $\vee$. Similarly, in the case of the resolution rule the operation performed on trees is deleting the branches of $p$ and $\neg p$ while keeping the rest intact.

% THIS MAY ONLY CONFUSE THE READER:

%% Next we consider an application of cloning of the form
%% \[
%% \frac{A\vee B\vee C\vee\dots}{(A\wedge A)\vee B\vee C\vee\dots}.
%% \]
%% Here $A$ must be exactly $\Pi^s_{k-1}$. If it were exactly $\Sigma^s_k$, the first formula would be fine, but the second would not be well-formed.

\paragraph{Negations.}\label{negations} In order to simulate the sequent calculus we need to define negations of stratified formulas. Given a $\Pi^s_k$ formula $A$, if we just dualize it, as we did before, we get a $\Sigma^s_k$. These formulas cannot occur as subformulas in the same proof, because the bottom connectives are different. Therefore we define the stratified negation of $A$ to be the dualized formula with literals padded by the bottom connective of $A$. So $\neg A$ is a stratified $\Sigma^s_{k+1}$ formula.

\begin{quote}
{\it Example} If $A$ is $p\wedge\neg q$, then $\neg A$ is $\wedge(\neg p)\vee\wedge(q)$.
\end{quote}

\paragraph{Efficient simulations.} Given two formula schemas $\cal A$ and $\cal B$, we will say that ${\cal A}\vdash{\cal B}$ can be \emph{efficiently simulated} if there exists a polynomial simulation which, for every given instance $A \vdash B$, produces a proof in bounded depth Symmetric Calculus in which the depth of the formulas does not exceed the depth of the two formulas, i.e., if $A,B\in\Pi^s_i$ (respectively $A,B\in\Sigma^s_i$), then all formulas in the proof are in $\Pi^s_i$ (in $\Sigma^s_i$).

\subsection{Zipping and unzipping}

There is ambiguity in representing formulas by stratified formulas, one formula may have several representations. We must show that it is easy to transform one representation to any other representation, otherwise the system would not be natural. The basic transformations that enable us to do this will be called zipping and unzipping.

\ex{Suppose we need to represent the formula $p\vee q$ as a $\Sigma^s_3$ formula. Then we have two possibilities: (1) $\vee\wedge(p\vee q)$, (2) $(\wedge\vee p)\vee(\wedge\vee q)$.}
{\small
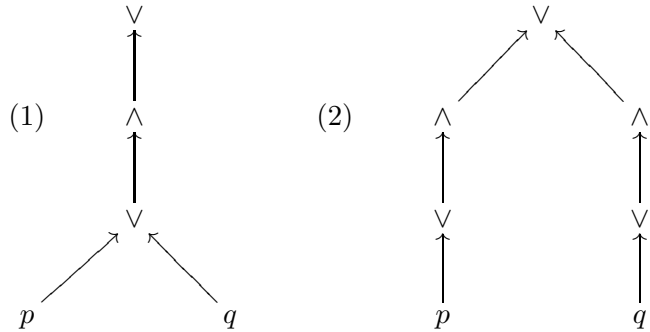
\begin{displaymath}
  \xymatrix{
& {\vee} & && &\vee&\\
(1)&\wedge\ar[u]& &(2)& \wedge\ar[ur]&&\wedge\ar[ul] \\
& \vee\ar[u] & && \vee\ar[u]&&\vee\ar[u] \\
p \ar[ur]& & q \ar[ul] && p \ar[u]& & q \ar[u]\\
  }% xymatrix
\end{displaymath}
}
% aby to fungovalo, nesmi se vypoustet radky!

\medskip
In general we can have arbitrary formulas instead of literals. 
We call the operation (1)$\mapsto$(2) \emph{unzipping} and the converse \emph{zipping}. 
We can zip and unzip more than two vertices, but it always has to be an even number, because the connective must be preserved.

\bfa
Any representation of a formula can be transformed to any other representation by zipping and unzipping. The number of the operation is bounded by the size of the formulas.
\efa
The proof is an easy exercise (which also involves a proper definition of representation).

\bll{l-zip}
Zipping and unzipping can be efficiently simulated, assuming that cuts and dual cuts can be.
\el
\bprf
Due to the symmetry of the calculus it is enough to simulate zipping and unzipping of conjunctions. Unzipping conjunctions is easy. It is done by cloning and weakenings:
\[
\wedge\vee(A\wedge B)\ \mapsto\ 
(\wedge\vee(A\wedge B))\wedge(\wedge\vee(A\wedge B))\ \mapsto\ 
(\vee\wedge(A))\wedge(\vee\wedge(B)).
\]
To simulate zipping we use a dual cut and a cut:
\[
(\vee\wedge(A))\wedge(\vee\wedge(B))  \mapsto\
((A\wedge B)\vee \neg B)\wedge(\vee\wedge(B))  \mapsto\
\wedge\vee(A\wedge B).
\]
% depth is OK
\eprf
%In this lemma we have assumed that cuts and dual cuts can be polynomially simulated within a given depth bound, which is yet to be proved.

Note that this implies that we can efficiently simulate zipping and unzipping not only conjunctions and disjunctions with two terms, but with any number greater or equal to~2.

%%%%%%%%%%%%%%%%%%%%%%%%%%%%%%%%%%%%%%%%%%%%%

\subsection{Cuts and dual cuts in the bounded depth system}

Now we want to show that one can efficiently simulate cuts and dual cuts in the bounded depth Symmetric Calculus. 
The simulation is the same as in the unbounded case except for one complication: zipped formulas. Suppose we want to simulate the cut
\[
\frac{(A\vee C)\wedge(\neg C\vee B)}{A\vee B}.
\]
If $C$ is a literal padded to the particular level, this is the resolution rule of the bounded depth Symmetric Calculus. If 
$C$ is $C_1\wedge C_2\wedge\dots\wedge C_n$ with $n\geq 2$, then we distribute $A$ to get $(A\wedge C_1)\vee(A\wedge C_2\wedge\dots\wedge C_n)$ and thus we recursively reduce the problem to simpler formulas $C_1$ and $C_2\wedge\dots\wedge C_n$. But if $C$ is a zipped conjunction $\wedge\vee\dots(C_1\wedge C_2\wedge\dots\wedge C_n)$, we cannot use distributivity immediately, we have to first unzip the formula. Unzipping conjunction is easy, but we also have to unzip the disjunction in $\vee\wedge\dots(\neg C_1\vee\neg C_2\dots\vee\neg C_n)\vee B$. To unzip disjunction, we need a dual cut and a cut. Fortunately, we need these operation for formulas that are simpler than $C$. So if we assume that we already have simulations for simpler formulas, we also have simulations of unzipping disjunctions for these formulas. 

So we prove the simulation by induction on the size of the cut and dual-cut formulas and get:
%logical depth of formulas. We have to use the \emph{true depth,} i.e., the depth of the formulas with the top unary conjunctions and disjunctions omitted. 

\bl
Cuts and dual cuts can be efficiently simulated in the bounded depth Symmetric Calculus.
\el

\subsection{Simulation of bounded depth sequent calculi}

Our aim now is to prove that our formalization of a bounded depth propositional proof is equivalent to the standard ones based on sequent calculi. In the Tait calculus restricted to depth $k$, sequents are sets of $\Pi_k$ and $\Sigma_k$ formulas. The interpretation of a sequents is the disjunction of the formulas in it, so they represent $\Sigma_{k+1}$ formulas. A set of sequents that appear in a proof can be represented by a conjunction of the $\Sigma_{k+1}$ formulas. Thus we get a formula of complexity at most~$\Pi_{k+2}$. Therefore the Symmetric Calculus system corresponding to a depth~$k$ sequent proof system is the system based on $\Pi^s_{k+2}$ proofs. 

\bprl{simulation}
For every $k\geq 1$, $\Pi^s_{k+2}$-Symmetric Calculus is polynomially equivalent to depth $k$ Sequent Calculus.
\epr
\bprf ~

1. We will show how to simulate the depth $k$ Tait Calculus. The Tait Calculus is formalized as follows (for more details, see~\cite{buss-intro}). Negations are only at variables, so $\neg A$ is the same formula as we defined in subsection~\ref{s-4.1}. The logical axioms are the sequents of the form $\Gamma,p,\neg p$, the rules are disjunction introduction, conjunction introduction, and cut. They will be simulated by dual resolution, weakening, the derived rule for distributive law of the form (\ref{e-distr}), and cut respectively. We have shown how to simulate the distributive laws and cut. It only remains to explain the technical issue concerning the representation of formulas and sequents. The complication is that we have to use stratified formulas.

For $\Pi_k,\Sigma_{k-1},\Pi_{k-2},\Sigma_{k-3},\dots$ formulas, we will use $\Pi^s_k,\Sigma^s_{k-1},\Pi^s_{k-2},\Sigma^s_{k-3},\dots$ stratified formulas; a particular representation is not important, because any such representation can be transformed to any other. For $\Sigma_k,\Pi_{k-1},\Sigma_{k-2},\Pi_{k-3},\dots$ formulas, we will use $\Sigma^s_{k+1},\Pi^s_{k},\Sigma^s_{k-1},$ $\Pi^s_{k-2},\dots$ stratified formulas obtained by shifting a stratified representation to higher level by adding unary conjunctions or disjunctions to literals. Negations $\neg A$ are defined as in section~\ref{negations}.

A sequent $\Gamma$ will be represented as follows. First we represent formulas of $\Gamma$ as described above. Then we pad every formula that has complexity smaller than $\Pi_k$ to level $\Pi_k$ by adding unary $\wedge$s and $\vee$s, say, on the top. Finally, we form a $\Sigma^s_{k+1}$ disjunction from these formulas. If the sequent consists of a single $\Pi_k$ formula, or formula of smaller complexity, the top disjunction in the $\Sigma^s_{k+1}$ stratified formula will be unary. The $\Sigma^s_{k+1}$ formulas representing $\Sigma_{k}$ formulas are disjunctions of $\Pi^s_k$ formulas and as such they will become parts of the  $\Sigma^s_{k+1}$ disjunction representing the sequent.

The simulation of a sequent proof is as follows. We represent the initial sequents by a conjunction of the stratified $\Sigma^s_{k+1}$ formulas. 
Then for every line in the proof we take the conjunction of all stratified formulas that represent sequents derived up to this line. If needed, we insert between two consecutive formulas a proof that simulates conjunction introduction, or cut. Finally, we use weakenings to remove all formulas except the one that we want to prove. It will be, of course, padded to the level~$\Pi^s_{k+2}$.

\medskip
2. Now we consider the opposite simulation. Given a $\Pi^s_{k+2}$ proof $A_1\dts A_n$, we transform every formula $A_i$ of the proof into a set of sequents $S_i$ with formulas of complexity at most $\Pi_k$: we interpret the $\Sigma^s_{k+1}$ disjunctions of $A_i$ as sequents of formulas obtained from the $\Pi^s_k$ subformulas of these disjunctions by omitting the unary $\wedge$s and $\vee$s if there are any. The disjunctions may contain multiple copies of $\Pi^s_k$ formulas, but this will not be reflected in the sets~$S_j$; they are sets, not multisets.
Furthermore, we omit the truth constants $\bot$ and $\top$.
We will show that the sequents of $S_{j+1}$ can be proved from the sequents of $S_j$ by polynomial size proofs of depth $k$. 

First we consider the case when a rule $R$ of the Symmetric Calculus is applied to a formula $A_i$ of the proof, not to its subformula. Since $A_i$ is a premise of a rule and its main connective is $\wedge$, $R$ can only be a rule from the right column. If $R$ is an instance of cloning there is no extra step in the simulation, because the sequent $S_j$ representing $A_i$ can be used repeatedly. We also do not have to simulate the rules for truth constants, because they do not appear in the sequents $S_j$. The $\wedge$-version of weakening (conjunction elimination) need not be simulated, because we keep all derived sequents in the Sequent Calculus.
%%  Dual resolution can only be applied in the form
%% \[
%% \frac{\top}{p\vee\neg p},
%% \]
%% where $\top$ is padded to a $\Pi^s_{k+2}$ formula, because if the premise were a conjunction $A\wedge\bot$ or $A\wedge B\wedge\bot$ with $A,B$ having complexity exactly $\Sigma^s_{k+1}$, the complexity of the  conclusion would exceed $\Pi^s_{k+2}$. So in this case we simulate the application of the rule by the axiom $p\vee\neg p$.
Finally, resolution is simulated by cut.

Suppose a rule $R$ of the Symmetric Calculus is applied to a $\Sigma^s_{k+1}$ subformula $B$ of $A_i$. Since the main connective of $B$ is $\vee$, $R$ can only be a rule from the left column. Contraction is simulated by contraction in the Sequent Calculus, $\bot$ is not used in the Sequent Calculus, and weakening is simulated by weakening in the Sequent Calculus (which is not among the rules of the Tait Calculus, but can be easily simulated).
The only rule that needs special treatment is dual resolution
\[
\frac{\bigwedge_i C_i\wedge\top\wedge\bigwedge_jD_j}{(\bigwedge_iC_i\wedge p)\vee(\bigwedge_jD_j\wedge\neg p)}.
\]
The formulas $\bigwedge_i C_i$ and $\bigwedge_jD_j$ are $\Pi^s_{k}$, hence they are simulated by formulas $C$ and $D$ of complexity at most $\Pi_k$. We need to show the following derivation:
\[
\frac{\Gamma,C\wedge D}{\Gamma,C\wedge p,D\wedge\neg p}.
\]
This is easy---it suffices to derive the sequent 
\[
\neg C,\neg D,C\wedge p,D\wedge\neg p,
\]
from which we get $\neg(C\wedge D),C\wedge p,D\wedge \neg p$ by disjunction introduction and then we can apply cut to get what we need: $C\wedge p,D\wedge\neg p$. To derive the sequent, first derive sequents $\neg C,C$ and $\neg D,D$ and $p,\neg p$. Then apply conjunction introductions.

When a rule of the Symmetric Calculus is applied to a subformula of complexity $\Pi^s_k$ or lower, then we proceed as follows. Let $B$ be the $\Pi^s_k$ subformula to whose subformula the rule is applied and let $C$ be $B$ after the rule is applied. Let $\hat B$ and $\hat C$ be the formulas that represent  $C$ and $B$ in the Sequent Calculus. Then we first prove an auxiliary sequent $\neg\hat B, \hat C$ and then we use cut to obtain $\hat C$.

We leave the construction of $\neg\hat B, \hat C$ to the reader. 
To construct $\neg\hat B, \hat C$ in the case when the rule is applied $B$ itself, not to a proper subformula of $B$, use the same argument as we used for dual resolution. If the application is deeper, use induction on the complexity of the formulas.
\eprf

\subsection{Deep inferences {using axioms}}

%Since we want to prove a precise characterization of interpolation pairs, we need precise bounds on the depth and size of some basic derivations. 

 Suppose we want to derive a formula from axioms. The axioms are typically clauses of a CNF formula from which we want to derive contradiction. So we start with the conjunction of the axioms and  gradually extend the conjunction by adding new derived formulas. In this process we need to do deep inferences \emph{using axioms}. E.g., we have an axiom $\neg A\vee B$ (representing $A\to B$) and there is an occurrence of $A$ deep inside of the currently derived formula. We want to replace $A$ by $B$ there. To this end we need to insert a copy of  $\neg A\vee B$ next to the occurrence of $A$, so that on the position of this occurrence we will get $(\neg A\vee B)\wedge A$. Then we can use cut to reduce it to $B$. The following lemma shows that a padded formula $A$ can be inserted into a disjunction.

\bll{l-insert}
Let $\vee(A)\wedge (B_1\vee B_2\vee\dots\vee B_n)$ and $(A\wedge B_1)\vee B_2\vee\dots\vee B_n$ be $\Pi^s_k$ formulas. Then there is a polynomial size depth $\Pi^s_k$ derivation of
\bel{e-ins}
\vee(A)\wedge(B_1\vee B_2\vee\dots\vee B_n)\ \vdash \ ((A\wedge B_1)\vee B_2\vee\dots\vee B_n).
\ee
\el
\bprf
This is just a weaker version of the distributive law, see Lemma~\ref{l-4.2}, and it is proved in the same way.
\eprf

By iterating this lemma, we may insert $A$ as deeper and deeper until all unary padding is removed from it.

%% We can easily simulate it using the distributive law, but then we cannot guarantee that we do not exceed the limit on the depth. So we will look more closely on the issue.
%% \bprf
%% The proof is essentially the same as for the distributive law except that we only apply it to formula $B_1$.
%% \[\begin{array}{ll}
%% A\wedge(B_1\vee B_2\vee\dots\vee B_n)&\\
%% A\wedge(\neg A\vee(A\wedge B_1)\vee B_2\vee\dots\vee B_n)&\mbox{ by dual cut}\\
%% ((A\wedge B_1)\vee B_2\vee\dots\vee B_n)&\mbox{ by cut.}
%% \end{array}\]
%% We only need to check that the middle formula is $\Pi^s_k$. Since we assume that the conclusion is $\Pi^s_k$, it must, in fact be $\Sigma^s_{k-1}$ and the true depth of $A$ can be at most $\Pi^s_{k-2}$. Thus a more precise representation of (\ref{e-ins}) is
%% \[
%% \vee(A)\wedge(B_1\vee B_2\vee\dots\vee B_n)\ \vdash \ ((A\wedge B_1)\vee B_2\vee\dots\vee B_n).
%% \]
%Hence $\neg A$ is $\Sigma^s_{k-2}$, which proves that the middle formula does not exceed $\Pi^s_{k-2}$.
%\eprf

%% \paragraph{Simulated substitution.}

%% \bll{l-subst}
%% Let $(\bigvee_i A_i)\wedge(\bigwedge_i\bigvee_j B_{ij})$ and $\bigvee_{ij}(A_{i}\wedge B_{ij})$ be $\Pi^s_k$ formulas. Then there is a polynomial size depht $\Pi^s_k$ derivation of
%% \[
%% (\bigvee_i A_i)\wedge(\bigwedge_i\bigvee_j B_{ij})\ 
%% \vdash\ \bigvee_{ij}(A_{i}\wedge B_{ij}).
%% \]
%% \el

%%%%%%%%%%%%%%%%%%%%%%%%%%%%%%%%%%%%%%%%%%%%%%%%%

\section{Games}\label{s4}

In this section we will introduce a new kind of games that we will use to characterize the interpolation pairs of bounded depth sequent calculi. First we describe the games in an intuitive way. The formal definition is in the next section.

We start with a concept that is a general form of many combinatorial games and it will be the bottom layer in our hierarchy of games. Such a game has two numerical parameters $n$, the length of the game and $m$, the number of symbols. In general, $m$ can be exponential in $n$, but we prefer to imagine that it is polynomially bounded. The actual relationship will depend on applications. A game of this type is given by
\ben
\item sets of symbols $A_1,A_2\dts A_n\sub[m]$,
\item transition functions $T_0:\{0,1\}\to A_1$, $T_i:\{0,1\}\times A_i\to A_{i+1}$, $i=1\dts n-1$, and
\item $W\sub A_n$, a set of of winning symbols.
\een
The game is played by two players, Player I and Player II, who alternate in choosing one of the the actions $T_0(0)$ or $T_0(1)$ and, for $i>1$, $T_i(0,x)$ or $T_i(1,x)$, where $x$ is the symbol that the previous player played. They play until they produce a sequence of symbols of length $n$. Player I wins if the last symbol played is an element of $W$, otherwise Player II wins.

It is possible to determine who has a winning strategy by a computation that runs in time polynomial in $n$ and $m$: for $i=n,n-1\dots 1$ compute inductively the set of winning symbols in step $i$. In fact, if we remove $W$ from such a structure, we can view it as a monotone Boolean circuit where the possible inputs are (strings encoding) sets $W\sub A_n$. For a given $W$, the circuit outputs $1$ iff in the game with $W$, Player I has a winning strategy.

In order to motivate our generalization we give an alternative definition of these games. Such a game will simply be given by a \emph{nondeterministic automaton}  $T$ with the set of states $A$ and a set of accepting states $W\sub A$. We will again assume that there are always only two possible action the automaton can do.  The players alternate in choosing one of the possible actions at each step. {This is not literally equivalent to the previous definition, because it corresponds to the case when $A_1=\dots=A_n$, but it is only a minor modification.}

Our games have another numerical parameter $k$, called \emph{the depth of the game}. The game will again be given by a nondeterministic automaton $T$, but now it will also use a tape with $n$ squares. $T$ reads the symbol on the currently visited square, moves to an adjacent square, reads the symbol there, and rewrites it. We assume that the only thing that the automaton can remember is the symbol that it has just read. 
$T$ starts at the leftmost square and moves to the right until it reaches the $n$th square. It ends there if $k=1$, otherwise it  reverses the direction and goes to the left. At the first square it stops  if $k=2$, otherwise it reverses the direction, and continues in this manner until it passes the tape $k$-times. As in the previous definition, players  alternate in controlling the automaton and the set of winning symbols is some subset $W$ of symbols.

%In this description it seems pointless for the automaton to read what it has just printed, but we want to formalize the game in such a way that all information is on the tape, which means that the automaton does not have memory.

In order to simplify the formal definition of the games, we will assume that the set of symbols are the same for each step, and we will have only two transitions functions, one for the directions from the left to the right, and one for the opposite direction, and the first symbol played will be fixed.

Clearly, the tape plays no role when $k=1$, but it is important if $k\geq 2$. Specifically, it is not possible to use the simple backtracking method to decide who has a winning strategy in polynomial time for $k\geq 2$. We think that, in fact, it is not possible to decide it using any polynomial time algorithm.

We will be interested in a special kind of strategies called \emph{positional strategies}. A positional strategy is a set of rules that instructs a player which action to choose based solely on the current state of the automaton and the content of the square it reads. This restricts the class of strategies significantly, so it is possible that a player has a general winning strategy, but no positional winning strategy. The advantage of positional strategies is that they have concise descriptions, polynomial in $n$ and $m$, where the degree depends on the depth $k$, and that, given such a strategy, one can check in polynomial time that it is a winning strategy. 

%Before we prove it, we give a formal definition of the games.

\subsection{Definition of games and positional strategies}\label{sec-def-games}

In order to simplify the formalization, we will assume, w.l.o.g., that the game starts by Player I rewriting the \emph{second} symbol on the tape.

\begin{definition}\label{d1}
A \emph{game of depth} $k$ is given by the length of a round $n$, the number of rounds $k$, a finite alphabet $A$ with a distinguished symbol $\Lambda\in A$, two functions, the \emph{transition functions} or \emph{legal moves} of the game,

\medskip $\overrightarrow{T}:\{0,1\}\times A\times A\to A$,
 
\medskip $\overleftarrow{T}:\{0,1\}\times A\times A\to A$, 

\medskip\noindent and a set of winning symbols $W\sub A$.
\edf
We will denote by $T$ the pair $(\overrightarrow{T},\overleftarrow{T})$.

\paragraph{The play---playing the game.}
The play starts with a string $\bar a$ of $\Lambda$s of length $n$, viewed as a tape with $n$ squares in which $\Lambda$ is printed. 
Player I starts by choosing $h\in\{0,1\}$ and replacing the second symbol by  $\overrightarrow{T}(h,\Lambda,\Lambda)$. Then players alternate and replace the $i+1$-st symbol $a_{i+1}$ of $\bar a$ (which is just $\Lambda$ in the first round) with $\overrightarrow{T}(h,a_i,a_{i+1})$, where $h$ is chosen by Player I if $i$ is odd and Player II if $i$ is even. They go on until the end of the string. If $k\geq 2$ they go back using $\overleftarrow{T}$ and so on. 
When they reverse direction, the last symbol, respectively the first symbol, is not rewritten again. This means that if they arrive at the end of the tape and the symbols on the tape are $a_1\dots a_{n-1}a_n$, then the player whose turn it is rewrites $a_{n-1}$ to  $\overleftarrow{T}(h,a_{n-1},a_n)$ for some $h\in\{0,1\}$ and they continue in the direction to the left. The same happens at the beginning of the tape.%
\footnote{This rule of the game is not essential, but it makes formalization simpler.}

After $k$ rounds the play ends and Player I wins if the last played symbol is in $W$, otherwise Player II wins.

\begin{definition}
A \emph{positional strategy for Player I} is a pair $(\overrightarrow{\sigma},\overleftarrow{\sigma})$ such that  

\medskip
$\overrightarrow{\sigma}:[k]\times [n]_{odd}\times A\times A\to A$,

\medskip
$\overleftarrow{\sigma}:[k]\times [n]_{odd}\times A\times A\to A$, 

\medskip
\noindent where $[n]_{odd}$ is the set of odd numbers $\leq n$. Furthermore, $\overrightarrow{\sigma}$ and $\overleftarrow{\sigma}$ must be compatible with~$T$, which means that for every $r,i,b,c$, there exists an $h$ such that $\overrightarrow{\sigma}(r,i,b,c)=\overrightarrow{T}(h,b,c)$ and similarly for $\overleftarrow{\sigma}$ and $\overleftarrow{T}$.

\medskip
A \emph{positional strategy for Player II} is a pair $\sigma=(\overrightarrow{\sigma},\overleftarrow{\sigma})$ defined in a similar way with $[n]_{odd}$ replaced by $[n]_{even}$.
\edf

For better readability, we will write $r$ and $i$ as subscripts, e.g., $\overrightarrow{\sigma}_{r,i}(b,c)$. The arrows above $\sigma$ are, clearly, determined by the index $r$, but we prefer to keep the arrows to stress the direction the strategy is used.

Instead of viewing the game as rewriting symbols on a tape, it is better to imagine that the players choose symbols in a $k\times n$ matrix in a zig-zag way, and the admissible choices are given by the previously played symbol and the symbol above the square that is to be filled. We will call such a partially filled matrix a \emph{history} of the play. More precisely, a \emph{history up to step $(r,i)$} is the record of a game played up to this step with the rest of the $k\times n$ matrix filled with $\Lambda$s.  We note a couple of useful properties of the history matrix $M$. 
\bel{e-Lambda}
M_{1,1}=\Lambda,
\ee
\bel{e-strana1}
M_{2r,n}=M_{2r-1,n}\mbox{ for all }1\leq r\leq k/2,
\ee
\bel{e-strana2}
M_{2r+1,1}=M_{2r,1}\mbox{ for all }1\leq r< k/2.
\ee
The latter two express that the players do not rewrite the last/first symbol when reversing the direction of playing. 

We will call a column vector of at most $k$ symbols of elements of $A$ a \emph{position} and view it as the first $r$, $r\leq k$, entries of a column of a $k\times n$ matrix. Given a history $M$ we say that $\bar a$ is the $(r,i)$-position if $\bar a$ is the column vector $(M_{1,i}\dts M_{r,i})^{\intercal}$.%
\footnote{$^\intercal$ denotes transposition of vectors (here a row vector to a column vector).} 

\ex{ Consider a 3-round game in which a play reached a position $(3,i)$. Let this be the history up to this position:

\[\baa{l|l|l|l|c|c|c|c|l|l|}
\cline{2-10}
&~\Lambda~  &M_{1,2}  &\dots&\rightarrow& \dots &{M_{1,i}} &{M_{1,i+1}} & \dots & M_{1,n} \\
\cline{2-10}
& M_{2,1}&\dots  &\dots&\leftarrow& \dots &{M_{2,i}} &{ M_{2,i+1}} & \dots & M_{2,n}(=M_{1,n})\\
\cline{2-10}
& M_{3,1}(= M_{2,1})&\dots &\dots&\rightarrow& \dots & {M_{3,i}} & & &\\
\cline{2-10}
\end{array}\]
Then the next symbol played, i.e.,  $M_{3,i+1}$, must be either $\overrightarrow{T}(0,M_{3,i},M_{2,i+1})$ or $\overrightarrow{T}(1,M_{3,i},M_{2,i+1})$. If moreover it is the turn of Player I and he uses a strategy $\sigma^I$, then $M_{3,i+1}=\overrightarrow{\sigma}_{3,i}^I(M_{3,i},M_{2,i+1})$ and the same for Player~II and strategy $\sigma^{II}$.
}
In general, $\overrightarrow{T}$ and $\overrightarrow{\sigma}_{r,i}$ are always applied to $M_{r,i},M_{r-1,i+1}$ (for $r$ odd), and $\overleftarrow{T}$ and $\overleftarrow{\sigma}_{r,i+1}$ are always applied to $M_{r-1,i},M_{r,i+1}$ (for $r$ even).

Let a game of depth $k$ be given. 
Let $r,s\leq k$, $i\leq n$, $\bar a\in A^r$, $\bar b\in A^s$. We say that $\bar a$ and $\bar b$ are 
$T$-\emph{compatible on steps $i$ and $i+1$}, if they are compatible with the transition functions applied at particular places. Given a strategy $\sigma$, we say that $\bar a$ and $\bar b$ are 
\emph{$\sigma$-compatible on steps $i$ and $i+1$}, if they are $T$-compatible and moreover they are compatible with $\sigma$ applied at particular places. 
A formal definition of $T$ and $\sigma$ compatibilities would be a long list of cases and formulas, while the concept is intuitively clear. Therefore we state formally only one case. 

If $\sigma$ is a strategy for Player I and $i<n$, $r\leq s$, $i$ and $r$ even, then  $\bar a$ and $\bar b$ are 
\emph{$\sigma$-compatible on steps $i$ and $i+1$} if there exits $h_1,h_3,\dots,h_{r-1}\in\{0,1\}$ such that 
\bel{e-comp}
\ba{lcl}
b_1&=&\overrightarrow{T}(h_1,a_1,\Lambda),\\
a_2&=&\overleftarrow{\sigma}_{2,i+1}(a_1,b_2),\\
b_3&=&\overrightarrow{T}(h_3,a_3,b_2),\\
a_4&=&\overleftarrow{\sigma}_{4,i+1}(a_3,b_4),\\
&\dots&\\
b_{r-1}&=&\overrightarrow{T}(h_{r-1},a_{r-1},b_{r-2}),\\
a_r&=&\overleftarrow{\sigma}_{r,i+1}(a_{r-1},b_r).
\ea
\ee
We will abbreviate {\it``$T$/$\sigma$-compatible on positions $i$ and $i+1$''} by {\it ``$T$/$\sigma,i$-compatible''} and omit $i$ if it is determined by the context.

For $r\leq k$, $i\leq n$, $\bar a\in A^r$ and $\sigma$ a strategy, we say that a position $\bar a$ is $(r,i,\sigma)$-\emph{reachable}, if there is a history of a play played according to $\sigma$ in which $\bar a$ is the $(r,i)$-position. We will denote by $R^\sigma_{r,i}$ the set of $(r,i,\sigma)$-{reachable} positions.

\bll{l-np}
For $k$ constant, 
given a positional strategy for Player I (respectively Player II), represented as a string of $kn|A|^2$ symbols from $A$, it is possible to decide in polynomial time if it is a winning strategy for Player I (respectively Player II).
\el
\bprf
We will show that sets $R^\sigma_{r,i}$ satisfy the following inductive conditions. 

For $r=1$,
\bel{e-c0}
{a}\in R^\sigma_{1,i+1}\ \equiv\ \exists {b}\in  R^\sigma_{1,i}({b},{a}\ \mbox{are  $\sigma$-compatible}).
\ee
If, e.g., $i$ is odd and $\sigma$ is strategy for Player I, this means that $\overrightarrow{\sigma}_{1,i}(a,\Lambda)=b$. 

For odd $r\geq 3$,
\bel{e-c1}
\bar{a}\in R^\sigma_{r,i+1}\ \equiv\ (a_1\dts a_{r-1})\in  R^\sigma_{r-1,i+1}
\wedge\exists \bar{b}\in  R^\sigma_{r,i}(\bar{b},\bar{a}\ \mbox{are $\sigma$-compatible}).
\ee

For even $r\geq 2$,
\bel{e-c2}
\bar{a}\in R^\sigma_{r,i}\ \equiv\ (a_1\dts a_{r-1})\in  R^\sigma_{r-1,i}
\wedge\exists \bar{b}\in  R^\sigma_{r,i+1}(\bar{a},\bar{b}\ \mbox{are $\sigma$-compatible}).
\ee

We will only prove (\ref{e-c1}), the proof of (\ref{e-c2}) is similar, and  (\ref{e-c0}) is trivial.
 
Suppose $\bar{a}\in R^\sigma_{r,i+1}$. Let $M$ be the history of a play in which $\bar{a}$ is the $(r,i+1)$-position. Let $\bar{b}$ be the $(r,i)$ position in $M$. Then, clearly, the right-hand side is satisfied.

We will use a ``hybrid argument'' to prove the opposite implication. Suppose that $(a_1\dts a_{r-1})\in R^\sigma_{r-1,i+1}$, $\bar{b}\in  R^\sigma_{r,i}$ and $\bar{b},\bar{a}$ are $\sigma$-compatible. We will consider histories of plays played according to $\sigma$. Let $M$ be a history of a play up to the point $(r,i)$ in which $\bar{b}$ is the $(r,i)$-position and $N$ be a history of a play up to the point $(r-1,n)$ in which $(a_1\dts a_{r-1})^{\intercal}$ is the $(r-1,i+1)$-position. Let $M'$ be the matrix consisting of the first $i$ columns of $M$ and $N'$ be the matrix consisting of the last $n-i$ columns of $N$. It is not difficult to see then that $M'N'$ is a history of a play up to $(r,i)$ in which $\bar{b}$ is the $(r,i)$-position and $(a_1\dts a_{r-1})^{\intercal}$ is the $(r-1,i+1)$-position. Now we can continue the play one more step to obtain $\bar a$ as the $(r,i+1)$-position. This finishes the proof of~(\ref{e-c1}).

Conditions (\ref{e-c1},\ref{e-c2}) give us a recursive procedure to compute the sets $R^\sigma_{r,i}$. The procedure has $kn$ steps and at each step we only need to consider at most $|A|^{2k}$ positions. Since $k$ is constant, this gives us a polynomial time algorithm. Finally we only need to check that the last set contains only winning positions of the player in question.
\eprf

\subsection{Modifications of the game}

There are various modifications of the definition of the game that are equivalent in the sense that they can efficiently simulate each other. 
First we note that we can assume w.l.o.g. that the transition functions also depend on the position on the tape. To simulate such a game by one whose transition function does not depend on the position, we let the players encode the position in the printed symbols. Thus if the game proceeds to the right and a player knows that he is on the $i$th position of the tape, he will encode the number $i+1$ in the symbol he will play.

We also do not have to insist that players alternate regularly. We may even allow steps that are done without players deciding anything. Furthermore, whose turn it is to move may also depend on the symbol to which they arrive.

Another modification, one of those that we are going to use, is that the players do not have to go to the ends of the tapes and can reverse the direction at other places on the tape.  If we want to satisfy the original definition, we may introduce an auxiliary symbol and let the players play this symbol until the end of the row and then back until they get to the place where they were supposed to pass to the next row.

We leave the formal statements and simulations to the readers, because they are easy, but may be complicated to write down formally.

\subsection{The disjoint NP pairs of the games and the main theorem}

We can now define the disjoint NP pairs of the games.

\bdf
For $k\geq 1$,
\[
A_k:=\{G\ |\ G\mbox{ is a game of depth $k$ in which Player I has a positional winning strategy}\},
\]\[
B_k:=\{G\ |\ G\mbox{ is a game of depth $k$ in which Player II has a positional winning strategy}\}.
\]
\edf
The disjointness is obvious, the membership in NP is the consequence of Lemma~\ref{l-np}.

For $k=1$, every strategy is positional and one can decide in polynomial time who has a winning strategy by backtracking winning positions. Also note that if we fix the transition function $\overrightarrow{T}$ and view sets $W$ as inputs, then such a game schema is essentially a monotone Boolean circuit that for a given $W$ decides who has a winning strategy (cf. Section~\ref{s8}).

For $k\geq 2$, one can easily construct games in which neither player has a positional wining strategy. For games of depth~2, it is open if one can decide in polynomial time who has a positional winning strategy given the promise that one of the players has such a strategy.

In the following two sections we will prove our main theorem.

\bt[Main Theorem]\label{t-main}
For $k\geq 1$, the pair $(A_k,B_k)$ is polynomially equivalent to the interpolation pair of the depth~$k-1$ sequent calculus.
\et
We have stated the theorem for the depth~$k-1$ sequent calculus, but in the proof we will use the $\Pi^s_{k+1}$-Symmetric Calculus. By Proposition~\ref{p2.1}, this theorem  also implies that for $k\geq 2$, $(A_k,B_k)$ is polynomially equivalent to the canonical pair of the depth~$k-2$ Sequent Calculus. 

To prove the theorem we need to show two reductions:
\ben
\item from   $(A_k,B_k)$ to the interpolation pair of the $\Pi^s_{k+1}$~Symmetric Calculus; this is Lemma~\ref{l5.1}, and
\item from the interpolation pair of $\Pi^s_{k+1}$~Symmetric Calculus to  $(A_k,B_k)$; this is Lemma~\ref{l6.1}.
\een

%%%%%%%%%%%%%%%%%%%%%%%%%%%%%%%%%%%%%%%%%%%%%%%%%%%%%%%%%%%%%%%%%%%%%%

\section{Proofs from games}\label{s5}

In this section we will construct, for every $k\geq 1$, a reduction from the pair $(A_k,B_k)$ to the interpolation pair of depth $\Pi^s_{k+1}$-Symmetric Calculus. 

\bll{l5.1}
Given a game $G$ of depth $k$, one can construct in polynomial time formulas $\Phi(\bar x)$ and $\Psi(\bar y)$ with disjoint sets of variables $\bar x$, $\bar y$, and a $\Pi^s_{k+1}$ refutation $D$ of $\Phi(\bar x)\wedge\Psi(\bar y)$ such that $\Phi(\bar x)$ is satisfiable when Player I has a positional winning strategy in $G$ and $\Psi(\bar y)$ is  satisfiable when Player II has a positional winning strategy in $G$. 
\el
Thus $(A_k,B_k)$ is reducible to the interpolation pair of $\Pi^s_{k+1}$-Symmetric Calculus and, by Proposition~\ref{simulation}, also to the interpolation pair of the depth $k-1$ sequent calculus.

We will prove this lemma by formalizing the statements that $\bar x$ is a positional strategy for Player~I and $\bar y$ is a positional strategy for Player~II, and constructing a depth $\Pi^s_{k+1}$ proof that it is impossible that both strategies are winning.

\subsection{The formula}

Let a game $G$ of length $n$ and depth $k$ be given. In order to formalize a strategy $\sigma$ we will use not only the strategy but also sets of reachable positions. Thus we will have variables both for elements of the strategies and elements of $R^I_{r,i}$ and  $R^{II}_{r,i}$, $r=1\dts k$, $i=1\dts n$.

We will use the convention that a propositional variable representing the truth of a relation $P(a_1\dts a_t)$ is denoted by $[P(a_1\dts a_t)]$ to represent propositions about strategies:  $[\sigma^I_{r,i}(a,b)=c]$, $[\sigma^{II}_{r,i}(a, b)=c]$. For  $\bar a\in R^I_{r,i}$ and  $\bar b\in R^{II}_{r,i}$, we will denote the variables simply by  $R^I_{r,i}(\bar a)$ and  $R^{II}_{r,i}(\bar b)$ in order not to overload notation with unnecessary symbols. One should keep in mind that in this notation the relation $R^I_{r,i}$ is indeterminate while elements $a_1\dts a_r$ are fixed, so the propositional variables are indexed by $r,i,a_1\dts a_t$; and this also concerns $R^{II}$, $\sigma^I$, and $\sigma^{II}$.
%%  and the sort of the propositional variables is given by the symbol $P$.

We will use $\times$ to refer either to Player I, or Player II and $*$ to refer to either  direction $\to$, or $\leftarrow$. For $\sigma$, we can omit the index referring to a player, because the player is determined by the index $i$ ($i$ odd is for I and $i$ even for II), and we can also omit arrows, because they are determined by the indices $r$ of rows.

We will use implication $A\to B$ to represent $\neg A\vee B$.

\paragraph{Variables of the formula}
\ben
\item $[{\sigma}_{r,i}(a,b)=c]$
for $r=1\dts k$, $i=1\dts n$, $a,b,c\in A$ (the variables for the strategies of Players I and II),
\item $R^I_{r,i}(\bar a)$, $R^{II}_{r,i}(\bar a)$ for $\bar a\in A^r$, $r=1\dts k$, $i=1\dts n$ (the variables for the sets of reachable positions).
\een

\paragraph{Clauses of the formula}
\ben
\item Clauses saying that ``$\sigma$ is a positional strategy''
\[
\bigvee_{c}[\sigma_{r,i}(a,b)=c]
\]
for $r=1\dts k$, $i=1\dts n$, and where  the disjunction is over $c$ such that $a,b,c$ are $T$ compatible.%
\footnote{The formula, in fact, expresses that $\sigma$s are \emph{total relations} defined properly, i.e., we do not formalize that they are \emph{functions}. Recall that ``$a,b,c$ are $T$ compatible'' means $T^*(h,a,b)=c$ for some $h\in\{0,1\}$.}

\item Clauses expressing (\ref{e-Lambda}), (\ref{e-strana1}), and (\ref{e-strana2}):
\bel{e-lambda}
R^I_{1,1}(\Lambda),\ R^{II}_{1,1}(\Lambda)
\ee
\bel{e-side1}
R^{\times}_{2r-1,n}(a_1\dts a_{2r-1})\equiv R^{\times}_{2r,n}(a_1\dts a_{2r-1},a_{2r-1}) \mbox{ for all }1\leq r\leq k/2,
\ee
\bel{e-side2}
R^{\times}_{2r,1}(a_1\dts a_{2r})\equiv R^{\times}_{2r+1,1}(a_1\dts a_{2r},a_{2r}) \mbox{ for all }1\leq r< k/2.
\ee

\item Clauses expressing the inductive conditions (\ref{e-c1}),(\ref{e-c2}) for $R^{\times}_{r,i}$.

First we need to express that $\bar a$ and $\bar b$ are $\sigma^{\times}$-compatible. Consider Player II, $i$ odd, $r$ even, and $\bar a,\bar b\in A^r$. Then ``$\bar a,\bar b$ are $\sigma^I$-compatible'' is defined by the conditions~(\ref{e-comp}). Hence if there are no $h_1,h_3,\dots,h_{r-1}\in\{0,1\}$ such that 
$b_1=\overrightarrow{T}(h_1,a_1,\Lambda)$,
$b_3=\overrightarrow{T}(h_3,a_3,b_2)$,
\dots,
$b_{r-1}=\overrightarrow{T}(h_{r-1},a_{r-1},b_{r-2})$,
then $a,b$ are not $\sigma^I$-compatible. Otherwise they are  $\sigma^I$-compatible iff  
$a_2=\overleftarrow{\sigma}^I_{2,i+1}(b_2,a_1)$,
$a_4=\overleftarrow{\sigma}^I_{4,i+1}(b_4,a_3)$,
\dots,
$a_r=\overleftarrow{\sigma}^I_{r,i+1}(b_r,a_{r-1})$.
Thus compatibility can be expressed by a conjunction of propositional variables
\[
[a_2=\overleftarrow{\sigma}^I_{2,i+1}(b_2,a_1)]\wedge
[a_4=\overleftarrow{\sigma}^I_{4,i+1}(b_4,a_3)]\wedge
\dots\wedge
[a_r=\overleftarrow{\sigma}^I_{r,i+1}(b_r,a_{r-1})].
\]
For Player II and other $r$ and $i$ it is similar.

This enables us to express formulas (\ref{e-c1}) and (\ref{e-c2}) by small propositional formulas, but not small CNFs. Fortunately, we only need implications from the right to the left. In the case of $r$ odd, i.e.~(\ref{e-c1}), it is 
\[
(a_1\dts a_{r-1})\in  R^\times_{r-1,i+1}
\wedge\exists \bar b\in  R^\times_{r,i}(\bar b,\bar a\ \mbox{are $\sigma,i$-compatible})\to \bar a\in R^\times_{r,i+1},
\]
which is equivalent to
\[
\forall \bar b((a_1\dts a_{r-1})\in  R^\times_{r-1,i+1}
\wedge \bar b\in  R^\times_{r,i}\wedge(\bar b,\bar a\ \mbox{are $\sigma,i$-compatible})\to \bar a\in R^\times_{r,i+1}).
\]
This can be represented by a set of clauses, one for every $b\in A^r$. The case of $r$ even is similar.

\ex{ Let $r$ and $i$ be odd, let $\bar a,\bar b\in A^3$. Suppose that $b_2=\overleftarrow{T}(h,a_2,b_1)$ for some $h\in\{0,1\}$. We consider the situation where $(b_1,b_2,b_3)^{\intercal}$ is the $i$th column of a history matrix and $(a_1,a_2,a_3)^{\intercal}$ is the $i+1$st column.
Then we have the following clause
\[
R^I_{3,i}(b_1b_2b_3)\wedge R^I_{2,i+1}(a_1a_2)\wedge
[a_1=\overrightarrow{\sigma}^I_{1,i}(b_1,\Lambda)]\wedge
[a_3=\overrightarrow{\sigma}^I_{3,i}(b_3,a_2)]\to
R^I_{3,i+1}(a_1a_2a_3).
\]
}

\item Clauses saying that the strategies of both players are winning. Strategy $\sigma^I$ is winning for Player I if all final positions reachable using $\sigma^I$ are winning. If $k$ is odd, this means
\[
\forall a_k\ ((a_1\dots a_k)\in R^I_{k,n}\to a_k\in W).
\]
This is expressed in the propositional calculus by
\bel{e-w1}
\bigwedge_{a_k\not\in W}\neg R^I_{k,n}(a_1\dts a_k)
\ee
Similarly, for Player II, still assuming $k$ odd, we get
\bel{e-w2}
\bigwedge_{a_k\in W}\neg R^{II}_{k,n}(a_1\dts a_k)
\ee
In the case of $k$ even, the index $n$ at $R$ is replaced by $1$.

\een

This defines CNF formulas $\Phi_k(\bar x)$ and $\Psi_k(\bar y)$, where $X$ are propositional variables 
$[\sigma^I_{r,i}(a,b)=c]$, $R^I_{r,i}(\bar a)$, and $Y$ are propositional variables 
$[\sigma^{II}_{r,i}(a,b)=c]$, $R^{II}_{r,i}(\bar a)$. The formula $\Phi_k(\bar x)\wedge\Psi_k(\bar y)$ expresses a contradictory fact that both players have positional winning strategies.

\subsection{The refutation of the formula}

We will now construct a $\Pi^s_{k+1}$ derivation of contradiction from the formula $\Phi_k(\bar x)\wedge\Psi_k(\bar y)$ defined above. 

\bll{l5.2}
One can construct in polynomial time a $\Pi^s_{k+1}$-refutation of the formula $\Phi_k(\bar x)\wedge\Psi_k(\bar y)$.
\el

This lemma implies Lemma~\ref{l5.1} because $\Phi_k(\bar x)$ (respectively $\Psi_k(\bar y)$) is a formalization the fact that Player I (Player II) has a winning strategy. 

Before going into details, we will explain the essence of the proof. Formula $\Phi_k(\bar x)\wedge\Psi_k(\bar y)$ says that it is impossible that both players have positional winning strategies. We use positional strategies because we need formulas of certain complexity, but, clearly, there cannot be \emph{any pair} of winning strategies for opposing players. The standard argument is that if we run the two strategies, then at the end only one player wins, so the two strategies cannot be winning. Let us try to formalize it and see why this argument cannot be used for depth~$d$ games for $d>1$. 

Let $S_{r,i}(x_1\dts x_r)$ denote that position $(x_1\dts x_r)^\intercal$ can be reached by playing strategies $\sigma^I$ and $\sigma^{II}$. To get a contradiction, we need to show that there exists a position $(x_1\dts x_d)^\intercal$ such that  $S_{d,n}(x_1\dts x_d)$ holds true if $d$ is odd, and $S_{d,1}(x_1\dts x_d)$ if $d$ is even. Clearly, we have $S_{1,1}(\Lambda)$, thus $\exists x\ S_{1,1}(x)$. Applying $\sigma^I$ we get $\exists x\ S_{1,2}(x)$, then using $\sigma^{II}$ we get $\exists x\ S_{1,3}(x)$ and so on until we obtain $\exists x\ S_{1,n}(x)$. If the depth of the game $d=1$, we are done, because a position in the last step cannot be winning for both players.

If $d>1$ we would like to continue. By definition, if $S_{1,n}(x)$, then $S_{2,n}(x,x)$, thus we have $\exists x,y\ S_{2,n}(x,y)$. Now we want to prove $\exists x',y'\ S_{2,n-1}(x',y')$. Let $x$ and $y$ be such that $S_{2,n}(x,y)$. We know that $\exists x'\ S_{1,n-1}(x')$, but this is not enough; we need an $x'$ such that $\sigma^\times(x')=x$ (where $x$ is $I$ or $II$ depending on whether $n$ is even or odd). So in order to be able to go back to a position in the first column, we need that
\[
\exists x_1\dots\exists x_n(S_{1,1}(x_1)\wedge x_2=\sigma^I(x_1)\wedge 
S_{1,2}(x_2)\wedge x_3=\sigma^{II}(x_2)\wedge\dots) .
\]
But if expressed as a propositional formula, it has exponential size, because the range of quantification is of size $|A|^n$.

What we can do instead is this. Observe that we have 
\bel{e.1.n}
\forall x\exists y( S_{1,n}(x)\to S_{2,n}(x,y)),
\ee 
because for a given $x$, we can take $y=x$. Suppose w.l.o.g.\- that $n$ is even. We can go back, to the left, with this formula. Suppose $S_{1,n-1}(x')$. Then $S_{1,n}(\sigma^\times(x'))$. From (\ref{e.1.n}), we get a $y$ such that  $S_{2,n}(\sigma^\times(x'),y)$ and then we can conclude $S_{2,n-1}(x',\sigma^{\times'}(y))$ using~(\ref{e-c1}) or (\ref{e-c2}). Thus we have shown $\forall x'\exists y'( S_{1,n-1}(x')\to S_{2,n-1}(x',y'))$. Repeating this argument we eventually get
\be
\forall x''\exists y''( S_{1,1}(x'')\to S_{2,1}(x'',y'')).
\ee
Now recall that we have $\exists x''\ S_{1,1}(x)$, so we get  
$\exists x''\exists y''S_{2,1}(x'',y'')$. 
If $d=2$, we get a contradiction, because $(x'',y'')^\intercal$ is the final position and as such it cannot be reached by both strategies. 

If $d>2$, this does not work, but we can use a formula with more quantifiers, specifically, formulas with $d$ alternating quantifiers. 

Let's have a look at the complexity of the formulas used in the proofs sketched above. For $d=1$, we are aiming at $\Pi^s_2$ proofs, but formulas $\exists x\ S_{1,i}(x)$ translate to $\Sigma^s_2$, because $S_{1,i}$ is $R^I_{1,i}\wedge R^{II}_{1,i}$. We cannot use that proof as it stands, but we can turn it around and argue contrapositively. We start with $\bigwedge_x\neg S_{1,n}(x)$ and proceed to the left. Thus we obtain $\bigwedge_x\neg S_{1,1}(x)$, from which we get contradiction using $S_{1,1}(\Lambda)$.

The proof for odd $d\geq 3$ is similar except that we have to use a formula with more alternations of $\wedge$s and $\vee$s; see formula~(\ref{nabla}) below.

For $d=2$, the proof above can be formalized as $\Pi^s_3$ proof. In general, for $d\geq 2$ even, we use formulas~(\ref{delta}). However, the complexity of these formulas does not guarantee that the proof has the same depth; it is necessary to check it, which is a little tedious, but not difficult.

The rest of this subsection is devoted to the proof of Lemma~\ref{l5.2}.

\bprf
Let a game of depth $k$ be given. First we observe that
from (\ref{e-w1}) and (\ref{e-w2})
 we get, by weakening, all clauses
\bel{e-con1}
\neg R^{I}_{k,n}(\bar a)\vee\neg R^{II}_{k,n}(\bar a)
\ee
for all $\bar a$ if $k$ is odd, and 
\bel{e-con2}
\neg  R^{I}_{k,1}(\bar a)\vee\neg  R^{II}_{k,1}(\bar a)
\ee
if $k$ is even. (We did not use these clauses to define our formula because they mix both types of variables.) We now consider cases according to the depth of the game.

\bigskip {\bf Case} $k=1$. We need a $\Pi^s_2$ refutation, which is essentially a Resolution refutation. We will omit the index $1$ of $R_{1,i}$ and the arrow  $\to$ above $\sigma$, because the direction does not change in this case; we will also abbreviate $\sigma(a,\Lambda)$ by $\sigma(a)$, because the second argument is always the same.

First we show by induction for $i=n,n-1\dts 1$ that clauses
\bel{e-re1}
\neg R^{I}_{i}(a)\vee\neg R^{II}_{i}(a)
\ee
are derivable for all $a\in A$.

We already have it for $i=n$ by (\ref{e-con1}).

Suppose we have (\ref{e-re1}) for $i+1$ and we want to get it for $i$. Suppose moreover that $i$ is odd. Then $a,b$ is $\sigma^{I},i$-compatible if  $\sigma^{I}_i(a)=b$.  Hence clauses of the formula that represent inductive conditions are 
\bel{e-re2}
\neg  R^{I}_{i}(b)\vee\neg[\sigma^{I}_i(b)=a]\vee  R^{I}_{i+1}(a),
\ee
\bel{e-re3}
\neg  R^{II}_{i}(b)\vee  R^{II}_{i+1}(a)
\ee
for all $T$-compatible pairs $b,a$. From $\neg R^{I}_{i+1}(a)\vee\neg  R^{II}_{i+1}(a)$, (\ref{e-re2}), and (\ref{e-re3}), we get by resolution
\bel{e-re4}
\neg  R^{I}_{i}(b)\vee\neg  R^{II}_{i}(b)\vee\neg[\sigma^{I}_i(b)=a]
\ee
We also have 
\bel{e-re5}
\bigvee_a[\sigma^I(i,b)=a]
\ee
for all $b\in A$ where the disjunction is over all $a$ such that $b,a$ is $T$-compatible (see clauses of the formula). From (\ref{e-re4}), and (\ref{e-re5}) we get $\neg R^{I}_{i}(b)\vee\neg R^{II}_{i}(b)$ by several applications of resolution. Since for every $b\in A$ there is an $a\in A$ such that $b,a$ are $T$-compatible, we get this formula for all $b\in A$.
For $i$ even, the proof is analogous.

For $i=1$, (\ref{e-re1}) gives, in particular, 
\[
\neg  R^{I}_{1}(\Lambda)\vee\neg  R^{II}_{1}(\Lambda).
\]
Since our formula contains clauses $R^{I}_{1}(\Lambda)$ and $R^{II}_{1}(\Lambda)$,
we get a contradiction.

\bigskip
{\bf Case} $k\geq 2$ {\bf even.}  We will abbreviate by
\[
S_{r,i}(x_1\dots x_r)\ :=\ R^I_{r,i}(x_1\dots x_r)\wedge R^{II}_{r,i}(x_1\dots x_r).
\]
This formula expresses that $(x_1\dots x_r)^{\intercal}$ is reachable using both strategies. It is a $\Pi^s_1$ formula (a conjunction of literals). Further, we introduce an abbreviation for compatibility:
\[
C_{r,i}(x_1\dots x_r,y_1\dots y_r)\ \equiv_{def}\ 
x_1\dots x_r\mbox{ and }y_1\dots y_r\mbox{ are both $\sigma^I$ and $\sigma^{II}$ $i$-compatible.}
\]
Also this formula is $\Pi^s_1$. We will also need $C^I_{r,i}$ and $C^{II}_{r,i}$ representing $\sigma^I$, respectively $\sigma^{II}$ compatibility. So $C_{r,i}\equiv C^I_{r,i}\wedge C^{II}_{r,i}$.

\bl
Let $1\leq i<n$, $1\leq r<k$, $\bar u,\bar v\in A^r$, Then the following formulas have polynomial size $\Pi^s_3$ proofs from the formula $\Phi_k(\bar x)\wedge\Psi_k(\bar y)$. 
\bel{621}
S_{1,i}(a)\to\bigvee_b (C_{1,i}(a,b)\wedge  S_{1,i+1}(b))
\ee
for every $a\in A$;
\bel{622}
(C_{r,i}(\bar u,\bar v)\wedge S_{r+1,i}(\bar ua)\wedge S_{r,i+1}(\bar v))\to
\bigvee_b (C_{r+1,i}(\bar ua,\bar vb)\wedge  S_{r+1,i+1}(\bar vb))
\ee
for $r$ even and every $\bar u,\bar v\in A^r$,  $a\in A$;

\bel{623}
(C_{r,i}(\bar u,\bar v)\wedge S_{r+1,i+1}(\bar vb)\wedge S_{r,i}(\bar u))\to
\bigvee_a (C_{r+1,i}(\bar ua,\bar vb)\wedge S_{r+1,i}(\bar ua)))
\ee
for $r$ odd and every $\bar u,\bar v\in A^r$, $b\in A.$
\el

\bprf
We will only prove (\ref{622}); the other two can be proved in the same way. Let $r$ be odd and assume w.l.o.g. that $i$ is also odd. Then our formula contains clauses (inductive conditions on $R$)
\[
C^{I}_{r,i}(\bar u,\bar v)\wedge R^{I}_{r+1,i}(\bar ua)\wedge
R^{I}_{r,i+1}(\bar v)\wedge
[b=\overrightarrow{\sigma}^{I}_{r,i}(a,v_r)]
\to R^{I}_{r+1,i+1}(\bar vb),
\]
\[
C^{II}_{r,i}(\bar u,\bar v)\wedge R^{II}_{r,i}(\bar ua)\wedge
R^{II}_{r,i+1}(\bar v)\wedge
\to R^{II}_{r,i}(\bar ua).
\]
Further, by definition
\[
C^{I}_{r,i}(\bar u,\bar v)\wedge 
[b=\overrightarrow{\sigma}^{I}_{r,i}(a,u_r)]
\to C^{I}_{r+1,i}(\bar ua,\bar vb).
\]
We also have 
\[
C^{II}_{r,i}(\bar u,\bar v)\wedge 
[b=\overrightarrow{\sigma}^{I}_{r,i}(a,u_r)]
\to C^{II}_{r+1,i}(\bar ua,\bar vb),
\]
because $[b=\overrightarrow{\sigma}^{I}_{r,i}(a,u_r)]$ ensures $T$ compatibility. 
Using conjunction introduction, see~(\ref{e-distr}), we get for every legal $a$ and $b$,
\[
C_{r,i}(\bar u,\bar v)\wedge S_{r+1,i}(\bar ua)\wedge
S_{r,i+1}(\bar v)\wedge
[b=\overrightarrow{\sigma}^{I}_{r,i}(a,u_r)]\to C_{r+1,i}(\bar ua,\bar vb)\wedge  S_{r+1,i+1}(\bar vb).
\]
By resolving with $\bigvee_b[b=\overrightarrow{\sigma}_{r,i}^{I}(a,u_r)]$ we get for every $a$,
\[
C_{r,i}(\bar u,\bar v)\wedge S_{r+1,i}(\bar ua)\wedge S_{r,i+1}(\bar v)
\to
\bigvee_b (C_{r+1,i}(\bar ua,\bar vb)\wedge  S_{r+1,i+1}(\bar vb)).
\]
\eprf

For $i=1\dts n$, we will denote by $\Delta_i$ the following formula
\bel{delta}
\bigwedge_{x_1}(S_{1,i}(x_1)\to
\bigvee_{x_2}(S_{2,i}(x_1x_2)\wedge
\bigwedge_{x_3}(S_{3,i}(x_1x_2x_3)\to \dots
\bigvee_{x_k}S_{k,i}(x_1\dots x_k)))).
\ee
Note that $\Delta_i$ is $\Pi^s_{k+1}$. Our plan is:
\ben
\item Prove $\Delta_n$.
\item Construct proofs of $\Delta_{i+1}\vdash\Delta_i$ for $i=n-1\dts 1$. Thus we get $\Delta_1$.
\item Then it would suffice to prove $\neg\Delta_1$, but the complexity of this formula is $\Sigma^s_{k+1}$ which is too much. Instead, one can derive
\[%bel{delta-minus}
\bigvee_{x_1}(S_{1,i}(x_1)\wedge
\bigwedge_{x_2}(S_{2,i}(x_1x_2)\to
\bigvee_{x_3}(S_{3,i}(x_1x_2x_3)\wedge \dots
\bigvee_{x_{k-1}}S_{k-1,i}(x_1\dots x_{k-1}))))
\]%\ee
and use it with $\Delta_1$ to derive
\bel{neco}
\bigvee_{x_1,x_2\dts x_k} S_{k,1}(x_1,x_2\dts x_k).
\ee
But it is easier to derive (\ref{neco}) from $\Delta_1$ only using clauses of the formula $\Phi(\bar x)\wedge\Psi(\bar y)$.

\item Finally we derive contradiction from (\ref{neco}) using cuts with formulas (\ref{e-con2}), which are $\neg S_{k,1}(\bar a)$ in the new notation.
\een
Now we present the proofs of these four steps.

\bigskip
1. We will prove $\Delta_n$. From clauses (\ref{e-side1}) (where we only need the implications from the left to the right) we get, using conjunction introduction,
\[
\bigwedge_{x_1}(S_{1,i}(x_1)\to
(S_{2,i}(x_1x_1)\wedge
\bigwedge_{x_3}(S_{3,i}(x_1x_1x_3)\to 
(S_{4,i}(x_1x_1x_3x_3)\wedge
\dots
\bigvee_{x_k}S_{k,i}(x_1x_1\dots x_{k-1}x_{k-1})))).
\]
Then $\Delta_n$ follows by weakening (but we can use this formula to derive $\Delta_{n-1}$ as well).

\bigskip
2. We will now prove $\Delta_{i+1}\vdash\Delta_i$. This is more complicated and we need to use some abbreviations:
\[
A_r:=S_{r,i}(x_1\dots x_r),\ B_r:=S_{r,i+1}(y_1\dots y_r),\ C_r=C_{r,i}(x_1\dots x_r,y_1\dots y_r).
\]
First we rewrite formulas (\ref{621}), (\ref{622}), and (\ref{623}) using the abbreviations.
%; for technical reasons, we will include one of the formulas in the antecendent also into the consequnet.
\bel{631}
A_1\to\bigvee_{y_1}(B_1\wedge C_1)
\ee
\bel{632}
A_{r+1}\wedge B_r\wedge C_r\to\bigvee_{y_{r+1}}(B_{r+1}\wedge C_{r+1})
\ee
\bel{633}
A_{r}\wedge B_{r+1}\wedge C_r\to\bigvee_{x_{r+1}}(A_{r+1}\wedge C_{r+1})
\ee
These formulas have $\Pi^s_3$ proofs, so we can use them, because we are constructing a $\Pi^s_{k+1}$ proof where $k\geq 2$.

The following formulas are the first steps of the derivation of $\Delta_i$ with  $\Delta_{i+1}$ being the first formula (a).

{\footnotesize

\[\ba{lccc}
(a)\ \bigwedge_{y_1}(B_1\to
&\bigvee_{y_2}(B_2\wedge
&\bigwedge_{y_3}(B_3\to
&\bigvee_{y_4}(B_4\wedge\ 
\dots\\
\\
(b)\ \bigwedge_{x_1}\bigwedge_{y_1}((B_1\wedge C_1)\to
&\bigvee_{y_2}(B_2\wedge C_1\wedge
&\bigwedge_{y_3}(B_3\to
&\bigvee_{y_4}(B_4\wedge\ 
\dots\\ 
\\
(c)\ \bigwedge_{x_1}(A_1\to
&\bigvee_{y_1}
\bigvee_{y_2}(B_2\wedge C_1 \wedge
&\bigwedge_{y_3}(B_3\to
&\bigvee_{y_4}(B_4\wedge\ 
\dots\\ 
\\
(d)\ \bigwedge_{x_1}(A_1\to
&\bigvee_{y_1}
\bigvee_{y_2}(A_1\wedge B_2\wedge C_1 \wedge
&\bigwedge_{y_3}(B_3\to
&\bigvee_{y_4}(B_4\wedge\ 
\dots\\ 
&%\bigvee_{y_r}B_r)))))\\
\\
(e)\ \bigwedge_{x_1}(A_1\to
&\bigvee_{y_1}
\bigvee_{y_2}\bigvee_{x_2}(A_2\wedge B_2\wedge C_2 \wedge
&\bigwedge_{y_3}(B_3\to
&\bigvee_{y_4}(B_4\wedge\ 
\dots\\ 
\\
(f)\ \bigwedge_{x_1}(A_1\to
&\bigvee_{y_1}
\bigvee_{y_2}\bigvee_{x_2}(A_2\wedge B_2\wedge C_2 \wedge
&\bigwedge_{y_3}((B_3\wedge C_3)
\to
&\bigvee_{y_4}(B_4\wedge C_3\ 
\dots\\ 
\\
(g)\ \bigwedge_{x_1}(A_1\to
&\bigvee_{y_1}
\bigvee_{y_2}\bigvee_{x_2}(A_2\wedge B_2\wedge C_2 \wedge
&\bigwedge_{y_3}\bigwedge_{x_3}((A_3\wedge B_2\wedge C_2)
\to
&\bigvee_{y_4}(B_4\wedge C_3\ 
\dots\\ 
\\
(h)\ \bigwedge_{x_1}(A_1\to
&\bigvee_{y_1}\bigvee_{y_2}\bigvee_{x_2}(A_2\wedge
&\bigwedge_{x_3}(A_3
\to
&\bigvee_{y_3}\bigvee_{y_4}(B_4\wedge C_3\ 
\dots\\ 
\\
~~~~\dots&\dots&\dots&\dots
\ea\]
}
We will describe how these formulas follow from previous ones.

\medskip{\it Proof of $(a)\vdash(b)$.} Add $C_1$ using (\ref{e-ext-con}).

\medskip{\it Proof of $(b)\vdash(c)$.} Using cuts with (\ref{631}).

\medskip{\it Proof of $(c)\vdash(d)$.} Add $A_1$ using (\ref{e-ext-con}) and weakening. 

\medskip{\it Proof of $(d)\vdash(e)$.} Since want to keep $B_2$, we first clone it and then apply cuts with (\ref{633}).

\medskip{\it Proof of $(e)\vdash(f)$.} Add $C_3$ in the same way as in (b).

\medskip{\it Proof of $(f)\vdash(g)$.} Cut $B_3\wedge C_3$ with (\ref{633}) for $r=2$.

\medskip{\it Proof of $(g)\vdash(h)$.} Cut $B_2\wedge C_2$.
 
\medskip\noindent
At the end we get
\[\ba{lccll}
~~~~\dots&\dots&\dots&\bigwedge_{x_{k-1}} (A_{k-1}\to 
&\bigvee_{y_{k-1}}\bigvee_{y_k} (B_k\wedge C_{k-1}))\dots )\\
\\
~~~~\dots&\dots&\dots&\bigwedge_{x_{k-1}} (A_{k-1}\to 
&\bigvee_{y_{k-1}}\bigvee_{y_k} (A_{k-1}\wedge B_k\wedge C_{k-1}))\dots )\\
\\
~~~~\dots&\dots&\dots&\bigwedge_{x_{k-1}} (A_{k-1}\to 
&\bigvee_{y_{k-1}}\bigvee_{y_k} (A_{k}\wedge C_{k}))\dots )
\ea\]
We do not need $C_k$ anymore, so we remove it by weakening (or, better, we use (\ref{633}) without it). What we get is $\Delta_k$ with additional disjunctions $\bigvee_{y_1}\dts\bigvee_{y_k}$. Since the formulas do not depend on $y_1\dts y_k$ anymore, the elements of these disjunctions are identical, hence we can get rid of the disjunctions by contractions.

\bigskip
3. We now prove  (\ref{neco}) from $\Delta_1$ and the clauses
\bel{clauses}
S_{1,1}(\Lambda),S_{2,1}(\Lambda x_2)\to S_{3,1}(\Lambda x_2 x_2),
S_{4,1}(\Lambda x_2x_2x_4)\to S_{5,1}(\Lambda x_2 x_2x_4x_4),\dots\ .
\ee
Here are the first steps of the proof.
\[\ba{llll}
\bigwedge_{x_1}(S_{1,1}(x_1)\to
&\bigvee_{x_2}(S_{2,1}(x_1x_2)\wedge
&\bigwedge_{x_3}(S_{3,1}(x_1x_2x_3)\to 
&\bigvee_{x_4}S_{4,1}(x_1x_2x_3 x_4)))\wedge\dots\\
\\

&\bigvee_{x_2}(S_{2,1}(\Lambda x_2)\wedge
&\bigwedge_{x_3}(S_{3,1}(\Lambda x_2x_3)\to 
&\bigvee_{x_4}S_{4,1}(\Lambda x_2x_3 x_4)))\wedge\dots\\
\\

&\bigvee_{x_2}(S_{3,1}(\Lambda x_2x_2)\wedge
&\bigwedge_{x_3}(S_{3,1}(\Lambda x_2x_3)\to 
&\bigvee_{x_4}S_{4,1}(\Lambda x_2x_3 x_4)))\wedge\dots\\
\\

&
&
&\bigvee_{x_2}\bigvee_{x_4}S_{4,1}(\Lambda x_2x_2 x_4)))\wedge\dots
\ea\]
Thus we get, in fact, a stronger formula
\[
\bigvee_{x_2}\bigvee_{x_4}\dots\bigvee_{x_{k-2}}\bigvee_{x_k}
S_{k,1}(\Lambda x_2x_2x_4x_4\dots x_{k-2}x_{k-2}x_k).
\]
from which 
we get contradiction using cuts with formulas (\ref{e-con2}).

\bigskip
{\bf Case} $k\geq 3$ {\bf odd.} For $i=1\dts n$, we will denote by $\nabla_i$ the following formula
\bel{nabla}
\bigwedge_{x_1}(S_{1,i}(x_1)\to
\bigvee_{x_2}(S_{2,i}(x_1x_2)\wedge
\bigwedge_{x_3}(S_{3,i}(x_1x_2x_3)\to \dots
\bigvee_{x_{k-1}}(S_{k-1,i}(x_1\dots x_{k-1})\wedge
\bigwedge_{x_k} \neg S_{k,i}(x_1\dots x_{k}))))).
\ee
Note that $\nabla_i$ is $\Pi^s_{k+1}$. The proof is similar to the proof for the case of $k$ even with a few modifications.

First, to derive $\nabla_n$ we also need to use clauses (\ref{e-con1}), which we now denote by $\neg S_{k,n}(x_1\dots x_{k})$.

We will now describe how to derive $\nabla_i$ from $\nabla_{i+1}$. Using our abbreviations and writing $\neg B_k$ as $B_k\to\bot$, formula $\nabla_{i+1}$ becomes
\[
\bigwedge_{y_1}(B_1\to
\bigvee_{y_2}(B_2\wedge
\bigwedge_{y_3}(B_3\to
\dots
\bigvee_{y_{k-1}}(B_{k-1}\wedge
\bigwedge_{y_k}(B_k\to\bot)\dots)))).
\]
This is like $\Delta_{i+1}$ for depth $k+1$ (which is even) the only difference being that $B_{k+1}$ is $\bot$. The same holds true for $\nabla_i$, so we can proceed in the same way as in the case of $k+1$. Having $\bot$ instead of $B_{k+1}$ makes our task even easier.

Finally, we derive contradiction by resolving $\nabla_1$ with clauses (\ref{clauses}). Here is how it goes for $k=3$.
\[\ba{llll}
\bigwedge_{x_1}(S_{1,1}(x_1)\to&\bigvee_{x_2}(S_{2,1}(x_1x_2)&\wedge
              &\bigwedge_{x_3}\neg S_{3,1}(x_1x_2x_3)))\\ \\
&\bigvee_{x_2}(S_{2,1}(\Lambda x_2)&\wedge&\bigwedge_{x_3}\neg S_{3,1}(\Lambda x_2x_3))\\ \\
&\bigvee_{x_2}(S_{3,1}(\Lambda x_2x_2)&\wedge&\bigwedge_{x_3}\neg S_{3,1}(\Lambda x_2x_3))\\ \\
&\bigvee_{x_2}(S_{3,1}(\Lambda x_2x_2)&\wedge&\neg S_{3,1}(\Lambda x_2x_2))\\ \\
&&\bot&
\ea\]
We leave the generalization for all odd $k\geq 3$ to the reader. 

This finishes the proof of Lemma~\ref{l5.2}.
\eprf

%%%%%%%%%%%%%%%%%%%%%%%%%%%%%%%%%%%%%%%%%%

\section{Games from proofs}\label{s6}

In this section we prove the opposite reduction, i.e., we will reduce the interpolation pair of $\Pi^s_{k}$-Symmetric Calculus to the pair $(A_k,B_k)$ of the depth $k-1$ games.

\bll{l6.1}
For every $k\geq 2$, given a refutation $D$ of a CNF formula $\Phi(\bar x)\wedge\Psi(\bar y)$ in the $\Pi^s_k$-Symmetric Calculus, where the  sets of variables $\bar x$ and $\bar y$ are disjoint, one can construct in polynomial time a game $G$ of depth $k-1$ such that if $\Phi(\bar x)$ is satisfiable then Player I has a positional winning strategy, and if $\Psi(\bar y)$ is satisfiable then Player II has a positional winning strategy. Moreover, the positional winning strategies can be constructed in polynomial time from the satisfying assignments.
\el

Let a refutation $D$ of a CNF formula $\Phi(\bar x)\wedge\Psi(\bar y)$ in the $\Pi^s_k$-Symmetric Calculus be given. We will assume that the refuted CNF is represented by a $\Pi^s_k$ formula as follows. If $k$ is even, then the bottom connectives of $\Pi^s_k$ are disjunctions. So in this case, we will simply be pad the CNF on the top, which is  schematically represented by
\[
\wedge\vee\dots\wedge(\dots\vee p\vee q\vee\dots).
\]
If $k$ is odd, the bottom connectives are conjunctions. So we will first pad literals to conjunctions and then we pad it on the top, which is schematically represented by
\[
\wedge\vee\dots\wedge(\dots\vee(\wedge(p))\vee(\wedge(q))\vee\dots).
\]
Let a $\Pi^s_k$ refutation $D:=(\Phi(\bar x)\wedge\Psi(\bar y)=\Gamma_1\dts\Gamma_m=\bot)$ be given. We will first define a game with $k$ rounds and $m$ steps in each round and then show that last round can be omitted so that we obtain a game of depth $k-1$.

%% The first round starts with $\bot$ padded to the level $\Sigma^s_{k-1}$ and then palyers continue playing $\Sigma^s_{k-1}$ subformulas of $\Gamma_2\dts\Gamma_m$. In the second roound they play $\Pi^s_{k-2}$ subformulas of $\Gamma_m\dts\Gamma_1$, and so on, until the last row which is  filled with literals occurring in $\Gamma_1\dts\Gamma_m$, or $\Gamma_m\dts\Gamma_m$. We require that a formula played on a row $i$ and column $j$, $i>1$, is a subformula of the formula played in the same column above it, i.e., in the row $i-1$ and column $j$. Furthermore, it must be ``logically connected'' with the previously played formula, which we will explain below. Note that, in particular, all formulas in column $j$ are subformulas of $\Gamma_j$.

The game starts at the last column that is associated with the last formula of the proof, which is $\bot$ padded to the level $\Pi^s_k$. The game starts with $\bot$ padded to level $\Sigma^s_{k-1}$ and then players proceed by Player's selecting maximal $\Sigma^s_{k-1}$ subformulas of $\Gamma_{m-1},\Gamma_{m-2},\dots$, so they select disjunctions from conjunctions. At some point they go to next row (if there is any) and change direction. When going to the right they select $\Pi^s_{k-2}$ conjunctions from the selected $\Sigma^s_{k-1}$ disjunctions. Then at some point they go to the next row (if there is any) and change direction again, and so on until the bottom row. 

Players cannot select an arbitrary subformulas, but only those that are in a certain sense ``logically connected''. Before we define the rules for selecting subformulas, we state the rules for changing directions.

\paragraph{Changing the direction.} Players go to the next row and change the direction
\ben
\item when they hit padded $\bot$ going to the right,
\item when they hit padded $\top$ going to the left,
\item when they get to the padded refuted CNF while going left.
\een
Note that when they hit $\bot$ going to the right, next time when going to the right they cannot get beyond it, because they hit it again, or hit another one before that. The same holds true for going to the left and hitting $\top$. In particular, if they hit $\top$ then they will not get to the initial CNF anymore.

There is a number of properties of the game we are defining that are not in accord with the formal definition given in Section~\ref{s4}, one of which is 
the possibility of going to the next row before the play reaches an end of the row. This was not allowed by the definition of the games in Subsection~\ref{sec-def-games}, but we have mentioned that it is possible to simulate such more general games. Another small and inessential discrepancy is that we defined games so that in the first and the last columns the symbols are not rewritten when starting in the opposite direction, which is not literally true in the case of the last column where one padding of $\bot$ is removed and the same concerns the first column where the refuted CNF is. The least important fact is that we start from the last column instead of the first one.

\paragraph{Legal moves.} As stated above the formula played must be a stratified subformula of either the formula in the proof (in the first round) or the formula played in the previous round of particular depth. Viewing formulas as trees, it must be a node connected to the previously played node of the tree. If the played subformula is not involved in an application of a deep inference rule, then the next subformula played is uniquely determined---it is the same formula on the corresponding position. In such a case the play proceeds without any action of the players. We will now define legal moves when a rule is applied to the subformula, or the subformula is a result of such an application.

Suppose the play proceeds from the left to the right. Then, for some $i$, the formulas played are $\Pi^s_i$ subformulas of $\Sigma^s_{i+1}$ formulas. Note that only the rules in the left column (see Section~\ref{s4.2}) change the structure of disjunctions, so we only need to consider them. 

In the following formulas we assume that $\bot,p,\neg p$ are padded to the appropriate level.

\ben
\item Permutation of a disjunction or conjunction. Players do not make any decisions; the play proceeds to the corresponding term of the disjunction or conjunction.
\item Contraction,
{\Large $\frac{\dots B\vee{A\vee A}\vee C\vee\dots}{\dots B\vee{A}\vee C\vee\dots}. $}

When any of the two occurrences of $A$ was played in the previous move, the next is the occurrence of $A$ in the conclusion of the rule. So also in this case players do not act. 

\item $\bot$ elimination, {\Large $\frac{\dots B\vee{\bot}\vee C\vee\dots}{\dots B\vee C\vee\dots}$.} 

If $\bot$ was played, the play cannot continue in the direction to the right. So the direction reverses and starts with $\bot$ with one padding removed.
\item Weakening,
{\Large $\frac{\dots B\vee C\vee\dots}{\dots B\vee {A}\vee C\vee\dots}.$}

All formulas from the premise are present in the conclusion, so the same formula is played as in the previous move.

\item Dual resolution, 
{\Large $\frac{\dots C\vee{(A\wedge{\top}\wedge B)} \vee D\vee\dots}{\dots C\vee{(A\wedge p)\vee(B\wedge\neg p)}\vee D\vee\dots}.$}

If $A\wedge{\top}\wedge B$ was the previous played subformula, then the legal moves are either $A\wedge p$ or $B\wedge\neg p$. Which of the two is played is decided by the player who owns the literal. When going back, if $p$ is chosen from $A\wedge p$, then $\top$ must be played in the next move (and the direction must be reversed).
\een

Now suppose the play goes from the right to the left. By symmetry, this is the same, except that now the game may reach the initial formula. Since in this direction the subformulas played are $\Sigma^s_i$ for some $i$, the play arrives either at some unary $\vee$ by which the initial formula is padded, or to a clause. If it is the padded formula, one padding is removed. If a clause is reached, then the player who owns it chooses a literal and they reverse the direction of play, i.e., if the clause is from variables $\bar x$ Player I chooses a literal, if it is from variables $\bar y$, Player II chooses a literal.

\paragraph{Termination of the game; winning positions.}

The game ends when players hit $\bot$, $\top$, or the first column when passing the bottom row.
On the bottom row they play a literal. The player whose literal hits $\bot$, $\top$, or the first column loses the game.

\ex{

\[\baa{|l|l|c|c|c|c|l|l|}
\cline{1-7}
\dots&\leftarrow&\dots&\dots&\dots&\dots&\vee\wedge\vee(\bot) \\
\cline{1-7}
\dots&\rightarrow&\dots&
A\wedge(r\vee p)\wedge(\neg p\vee s)\wedge B
&A\wedge(r\vee\bot\vee s)\wedge B
&\dots&\wedge\vee(\bot) \\
\cline{1-7}
\dots&\leftarrow&\dots&r\vee p&r\vee\bot\vee s&\dots&\vee(\bot)\\
\cline{1-7}
\dots&\rightarrow&\dots&p&\bot&\mbox{\tiny game ends}&\\
\cline{1-7}
\end{array}\]

In the last step of this play we interpret $\bot$ as logically connected with $p$, therefore $\bot$ follows after $p$. If this was not in the bottom row (in which case $p$ and $\bot$ would be padded), then the direction would be reversed and the play would go on.
}

\paragraph{The abridged game.}
It is clear that the bottom row is superfluous: once they get to this row, the same literal is played until the end of the game. Hence, we can omit this row and define the terminating positions to be the positions where they are supposed to go to the bottom row, and declare the position to be losing for the player who owns the literal to be played. In this way, from a $\Pi^s_k$ proof, we obtain a game of depth $k-1$.

When $k$ is even, literals of the initial formula are padded to conjunctions. Thus a literal from a disjunction may be chosen when going to the $k-1$st round, but the game is not decided yet. The literal that eventually hits $\bot$ may be different.

\subsection{A winning strategy from a satisfying assignment}

To prove Lemma~\ref{l6.1}, it suffices now to prove the following.

\bll{l6.2}
Given a satisfying assignment for $\Phi(\bar x)$, one can construct in polynomial time a positional winning strategy for Player I. 
The same holds for a satisfying assignment to $\Psi(\bar y)$ and Player II.
\el
\bprf
For this proof, it will be convenient to consider the full game, not the abridged version. 

Let $\bar a$ be a satisfying assignment for $\Phi(\bar x)$. Recall that the only decisions that Player I has to make occur when the rules of resolution and dual resolution are applied and when the players arrive at a clause $C(X)$ of the refuted CNF. The strategy is defined as follows:
\ben
\item The play proceeds to the left and arrives at a resolution step with a variable $x\in\bar x$. Then Player I chooses the disjunction in which the literal is \emph{falsified} by~$\bar a$.
\item The play proceeds to the right and arrives at a dual resolution step with a variable $x\in X$. Then Player I chooses the conjunction in which the literal is \emph{satisfied} by~$\bar a$.
\item The play proceeds to the left and arrives at a clause $C(\bar x)$  of the CNF. Then Player I chooses a literal from $C(\bar x)$ that is \emph{satisfied} by~$\bar a$ (the literal may be padded to a conjunction if the bottom layer of connectives are conjunctions).
\een

This is, clearly, a positional strategy. We will show that it is a winning strategy. We will consider two cases.

\medskip
1. \emph{The lowest level of connectives are disjunctions.} Then the players pass the bottom row, the level of literals, in the direction to the right. So the play stops when they hit $\bot$. We should show that the literal with which they hit $\bot$ cannot be from $\bar x\cup \neg\bar x$. Suppose, by way of contradiction, that this literal is $x\in \bar x\cup\neg \bar x$. 

If they started the bottom row at its beginning, which is the refuted CNF, then it means that Player I has chosen a literal from some clause $C(\bar x)$. Since he uses the strategy described above, the literal is satisfied by $\bar a$. On the other hand, in the resolution step that produced this occurrence of $\bot$, he chose the direction with the falsified literal. Since what he plays on a lower row is a subformula of a formula on a higher row, this is the same literal. (More precisely, on the higher row it is this literal padded to the appropriate level.) So this is not possible.

\ex{
Suppose in the previous example Player I used this strategy and $p\in\bar x$. Then he chose $r\vee p$ because $p$ was falsified by $\bar a$. When he started the bottom row, he should have picked a literal that is satisfied by $\bar a$, which is not $p$. So this situation cannot occur.
}

If they started at some occurrence of $\top$, then $x$ must be one of the literals to which $\top$ split. Again, Player I chose the literal that is satisfied by $\bar a$, but $x$ should be falsified by $\bar a$ because of the resolution step that produced $x$.

\medskip
2. \emph{The lowest level of connectives are conjunctions.} Then the players pass the bottom in the direction to the left and the play stops when they hit either $\top$, or the beginning of the row. Suppose that the literal with which they hit $\top$ or the beginning is $x\in \bar x\cup\neg\bar x$. 
Since they always start the bottom row from some occurrence of $\bot$ and go left, $x$ must be falsified by $\bar a$. There are two cases:
\bi
\item[(i)] If they hit $\top$, then we get a contradiction as above. 
\item[(ii)] Suppose they hit the beginning. This means that in the previous round they started from the beginning of the row. Then one player had to choose a literal (padded to a conjunction) from a clause of the CNF. Since $x\in \bar x\cup\neg\bar x$, it was chosen by Player I. But then, according to the rules of the strategy, it must be satisfied by $\bar a$. This is a contradiction again.
\ei
\eprf

\ex{Consider the unabridged version of the game. 
Let $k$ be odd. So the bottom connective is conjunction and the literals of the initial formula are padded to conjunctions. Suppose Player I plays the strategy based on a satisfying assignment for $\Phi(\bar x)$. Suppose they start the last but one round by Player I choosing a literal $x$ from a clause of $\Phi(\bar x)$; he chooses a satisfied literal. Now they proceed to the right. The single element conjunction $\wedge(x)$ may be enlarged as they go on. The rule that enables it is dual resolution. Literals from both $X$ and $Y$ can be added, but Player~I only picks those that are satisfied. Eventually all but one, say $p$, are removed by weakening and they hit $\bot$. The literal $p$ cannot be from $\bar x\cup\neg\bar x$, because if it were, then on some previous row  Player~I would decide where to go from $\bot$ and he would choose the unsatisfied literal from the two options. So $p$ belongs to Player~II. Then Player~II loses, because $p$ is repeated on the bottom level all the way to the end of the game. The game will end at the position where $p$ was introduced by dual resolution. There the game will hit $\top$.
}

%%%%%%%%%%%%%%%%%%%%%%%%%%%%%%%%%%%%%%%%%%%%%%%%%%%%%%%%%%%%%%%%%%%%%%%%

\section{A stronger result}\label{s7}

\paragraph{Game schemas.} We will call a \emph{game schema} a system of rules $S$ that defines legal moves and a set of end positions. What a schema does not specify is which positions are winning for which player. We will furthermore require the end positions to be labeled by $0$s, $1$s, and variables $z_1\dts z_n$. Given a schema $S(\bar z)$ and an assignment $\bar a:\{z_1\dts z_n\}\to\{0,1\}$, we obtain a game $S(\bar a)$ where winning positions of Player~I are the positions labeled by $1$ and those labeled by $0$ are winning for Player~II.

\medskip
Thus a game schema can be used to define a total monotone Boolean function $f:\{0,1\}^n\to\{0,1\}$ if we consider general strategies and a partial monotone function if we consider only positional strategies. The value of the function is $1$ (respectively $0$), if Player~I (Player~II) has a winning strategy. We will mostly be interested in positional strategies.

In the case of games that we have introduced it is very easy to define the corresponding concept of a game schema---it suffices to omit the set of winning symbols $W$ from Definition~\ref{d1} and suitably label the elements of the set $A$. We will call resulting objects \emph{depth $k$ game schemas}.

Our aim is to prove the following strengthening of polynomial simulation of interpolation pairs, which generalizes monotone feasible interpolation for Resolution.

%%%%%%%%%%%%%%%%%%%%%%%%%%%
\btl{t7.1}
Let $\Phi(\bar x,\bar z)$ and $\Psi(\bar y,\bar z)$ be two CNF formulas whose only common variables are $\bar z$. Suppose variables $\bar z$ occur in $\Phi$ only positively and in $\Psi$ only negatively.
Let a $\Pi^s_k$ refutation $D$ of  $\Phi(\bar x,\bar z)\wedge\Psi(\bar y,\bar z)$ be given, $k\geq 2$. Then it is possible to construct in polynomial time a depth $k-1$ game schema $S(\bar z)$ such that for every assignment $\bar a:\bar z\to\{0,1\}$, if  $\Phi(\bar x,\bar a)$ is satisfiable, then Player~I has a positional wining strategy in $S(\bar a)$ and  if  $\Psi(\bar y,\bar a)$ is satisfiable, then Player~II has a positional wining strategy in $S(\bar a)$.
\et
Note that in particular the size of the game schema $S$ is polynomial in the size of the proof $D$.
%%%%%%%%%%%%%%%%%%%%%%%%%%
\bprf
The proof is a simple adaptation of the proof of Lemma~\ref{l6.2}. Let a $\Pi^s_k$ refutation $D$ be given. We introduce new variables $z'_1\dts z'_n$ and substitute them for $\neg z_1\dts \neg z_n$ in $\Psi(\bar y,\bar z)$. Let $\Psi'(\bar y,\bar z')$ be the formula after the substitution. We take the CNF formula  $\Phi(\bar x,\bar z)\wedge\Psi'(\bar y,\bar z')\wedge\Delta(\bar z,\bar z')$, where  $\Delta(\bar z,\bar z')$ is the conjunction of all clauses $\neg z_i\vee\neg z'_i$. One can, clearly, construct a $\Pi^s_k$ refutation $D'$ of  $\Phi(\bar x,\bar z)\wedge\Psi'(\bar y,\bar z')\wedge\Delta(\bar z,\bar z')$ that is only slightly larger than $D$.

We define a game schema from $D'$ in the same way as we did in Lemma~\ref{l6.2} with one modification that concerns clauses $\neg z_i\vee\neg z'_i$. When the play arrives at such a clause, we let Player~I choose a literal from it. If he chooses $\neg z_i$ then the play continues in the usual manner. If he chooses $\neg z_i'$, then Player~II will have an opportunity to challenge Player~I's move. If she challenges, then the game ends and the end position gets label $z_i$. If she does not challenge, the game continues as before.

Let an assignment  $\bar a:\bar z\to\{0,1\}$ be given and suppose $\Phi(\bar x,\bar a)$ is satisfiable. Let $\bar b:\bar x\to\{0,1\}$ be the satisfying assignment. We will use $\bar a$ also for variables $\bar{z}'$ as if they were $\neg\bar z$. 
Player~I will use $\bar b,\bar a$ for his strategy in the same way as it was in Lemma~\ref{l6.2}. He controls the choice from clauses $\neg z_i\vee\neg z'_i$ and will always pick the satisfied literal. Player~II can challenge only if she picks $\neg z'_i$. But then $\neg z'_i$ is satisfied so $z'_i$ not satisfied and $z_i$ is satisfied, i.e., $a_i=1$. Thus this action of Player~II would result in Player~I immediately winning.

Suppose now that $\Psi(\bar y,\bar a)$ is satisfiable. This means that $\Psi'(\bar y,\neg\bar a)$ is satisfiable. 
Now Player~II does not control the action when they arrive at $\neg z_i\vee\neg z'_i$. But what she only needs for her strategy is that $\neg z'_i$ is not chosen if $\neg z'_i$ is not satisfied. If Player~I chooses $\neg z'_i$ in spite of $\neg z'_i$ being not satisfied, then $z'_i$ is satisfied and $z_i$ is not satisfied, i.e., $a_i=0$. Hence Player~II challenges and wins immediately. 
\eprf

Let $\cal S$ be a family of game schemas with a concept of a positional winning strategy. Then for $S(\bar z)\in{\cal S}$ we will denote by $f_S(\bar z)$ the partial Boolean function defined by $S$. Let $\cal S$ and $\cal T$ be two families of game schemas. Then we will say that $\cal S$ \emph{strongly polynomially} (respectively \emph{quasipolynomially}) reduces to $\cal T$, if for every game schema $S(\bar z)\in {\cal S}$, there exists at most polynomially (quasipolynomially) larger game schema $T(\bar z)\in{\cal T}$ such that $f_S\sub f_T$ (i.e., whenever $f_S(\bar z)$ is defined, so is $f_T(\bar z)$ and $f_S(\bar z)=f_T(\bar z)$). The word ``strongly'' refers to the fact that in the reductions the sets of variables are exactly the same.

\bco
Let $k\geq 2$. Suppose $\Pi^s_k$-Symmetric Calculus can quasipolynomially simulate $\Pi^s_{k+1}$-Symmetric Calculus on CNFs. Then game schemas of depth $k$ strongly reduce to game schemas of depth $k-1$.
\eco
We have stated the corollary for quasipolynomial simulation, because we know that the depth $k-2$ sequent calculus does not polynomially simulate the depth~$k-1$ sequent calculus (cf. \cite{impagliazzo-krajicek} and \cite{krajicek}, Theorem 14.5.1), hence also $\Pi_k$-Symmetric Calculus does not polynomially simulate $\Pi_{k+1}$-Symmetric Calculus.

%%%%%%%%%%%%%%%%%%%%%%%%%%%%%%%%%%%%%%%%%%%%%%%%%%%%%%%%%%%%%%%%%%%%%%%%%

\section{Two special cases}\label{s8}

We will consider two special cases: depth 1 and depth 2 games. We will show that depth 1 game schemas are essentially monotone Boolean circuits and depth~2 games are equivalent to point-line games introduced in~\cite{BPT}.

\subsection{Depth 1 games}

Let $C$ be a monotone Boolean circuit. $C$ is given by a directed oriented graph $H$ with a root~$r$. The root is the output of the circuit. Vertices are labeled by $\vee$ and $\wedge$, except for leaves which are labeled by variables $x_1\dts x_n$. Let $\bar a$ be an assignment to the variables. Then we can view the pair $C,\bar a$ as a game with two players $\bigvee$ and $\bigwedge$. They start at the root and follow the arrows with the direction chosen by the player by whose label the vertex is labeled. Player~$\bigvee$ wins iff they reach a leaf whose variable is substituted by $1$. It is not difficult to see that Player~$\bigvee$ has a winning strategy iff $C(\bar a)=1$. 

Thus in our terminology a monotone Boolean circuit is a game schema. It is a universal model in the following sense. If the number of configurations in a finite game schema $S$ is $N$, then it can be represented by a monotone Boolean circuit with $N$ vertices. In particular, monotone Boolean circuit can represent our depth 1 game schemas with the number of vertices polynomial in the size of the game schema. A corollary of this is:

\bprl{p9.1}
~

1. In every depth 1 game one of the players has a positional winning strategy.

2. One can decide in polynomial time who has a winning strategy in a depth 1 game.
\epr

The converse simulation is also easy. We only need to represent vertices of the graph of $C$ by elements of $A$, the set of symbols used in the game. Of course, in the definition of a circuit one does not require that $\vee$s and $\wedge$s alternate regularly, but this can easily be simulated by some dummy moves. Thus the power of depth 1 game schemas and monotone Boolean circuits is the same up to a polynomial increase.

\subsection{Depth 2 games}

The canonical NP pair of Resolution, which is polynomially equivalent to the interpolation pair of the depth 1 sequent calculus, has been characterized by a game called the \emph{point-line game} \cite{BPT}. It follows that the pairs of the point-line game and depth 2 game are polynomially equivalent. Here we will show direct simulations which also preserve monotonicity.

A point-line game is given by a directed acyclic graph $H$ with a root $R$ and some additional structure. We will view the nodes of the graph as having some  inner structure---like circles in which points are drawn. Each node is assigned either to player Black or player White. The root is empty, the other nodes contain some points and each leaf contains exactly one point. If there is an arrow from a node $P$ to a node $Q$, then there is a partial matching $M_{P,Q}$ between the points of $P$ and $Q$. A play starts at the root and proceeds along the arrows to a leaf. At each node the player who owns it decides where to proceed. When a node is visited, it is filled with black and white pebbles put on the points for the node. The configurations of pebbles are determined by the way in which the node was reached. The rule for pebbles is:
\bi
\item if the play goes $P\to Q$, then pebbles that are in the domain of  $M_{P,Q}$ are moved from $P$ along the lines to $Q$; the remaining points of $Q$ are filled with pebbles of the player other than the one who did this move.
\ei
When the play reaches a leaf, then the color of the pebble that ends up there decides who wins. 

\medskip
There are certainly many modifications that result in essentially the same concept. For us, the most important one is to allow more points in the root. This version is then a game schema where an instance is given by putting some pebbles on the points of the root. Then we can use this schema to compute partial monotone functions.

Once we allow points in the root, we can also w.l.o.g. assume that the range of each matching $M_{P,Q}$ covers all points in $Q$. Then the rule about pebbles becomes simpler---just move pebble along the matchings.

\medskip
Since both the depth 2 games and point-line games characterize the interpolation pair of the depth~1 sequent calculus, which is equivalent to $\Pi_2$-Symmetric Calculus, the NP-pairs of the two kinds of games are polynomially reducible to each other. Below we will show direct simulations that, moreover, show strong polynomial reduction between the corresponding game schemas.

\paragraph{Simulation of point-line games by depth 2 games.}
In order to see the connection with depth 2 games, consider a point-line game, with the modifications mentioned above, presented in a different way. The play starts with some configuration of pebbles on the root. Then the players traverse the graph, but instead of putting pebbles on the nodes they only mark the path they have taken. When they reach a leaf, they start on the point in that leaf and go back along the lines that connect points of the taken path. Thus they get back to a point in the root. The color of the pebble that is there decides who wins.

It is clear that the game is the same. Moving the pebbles in the original way of playing the game is only a means to save the trip back to the root.

In this formulation it is clear that the point-line game \emph{is} a a depth 2 game; no comment is needed. In fact, it is a special case of depth 2 game---in the last round players do not decide anything. This may suggest that the point-line game is a weaker concept, but this is not the case. We will show below that one can also simulate depth 2 games by point-line games.

\paragraph{Simulation of depth 2 games by point-line games.}
The idea of the simulation is to use the directed graph of the positions in the first round as nodes of the point-line game and positions in the second round as their points. A position in the second round in the $i$th column is a pair $(a,b)$ where $a$ is a symbol played in the first round and $b$ is a symbol played in the second round. Thus $a$ will be a node and $(a,b)$ a point in it. Black pebbles will represent winning positions of Player I and white pebbles the winning positions of Player II. This means that if $\Lambda,a_2\dts a_i$ are the first $i$ moves in the depth~2 game, then a black pebble on $(a_i,b)$ means that Player~I has a strategy to win the game if started from position $(a_i,b)$ and $\Lambda,a_2\dts a_i$ are fixed on the first row. Imagine that the game has been played until this point and it remains $i$ steps to finish the game. Having such winning positions for $(a_i,b)$ and $\Lambda,a_2\dts a_i$ and given $a_{i+1}$ that is a legal move after $a_i$, we can easily determine winning positions $(a_{i+1},b)$ for $\Lambda,a_2\dts a_i,a_{i+1}$. The process of defining the winning position $(a_{i+1},b)$ can be represented as moving black pebbles from node $a_i$ to node $a_{i+1}$, except that it is slightly more complicated than just moving pebbles along lines. Eventually we arrive at $a_n$. Then, by definition, there is only one position with $a_n$ as the first coordinate, namely $(a_n,a_n)$. If it is a winning position for $\Lambda,a_2\dts a_n$, which is represented by a black pebble on it, then Player~I has a strategy to win the game from this position with $\Lambda,a_2\dts a_n$ in the first row.

We will now describe the simulation in more detail. 
Let a depth 2 game be given. We will first construct a point line game with slightly more complicated rules for moving pebbles. We will allow  conjunctions and disjunctions with pebbles, which means that for two nodes connected by an arrow $P\to Q$ we may have two points $p_1,p_2\in P$ connected by lines to a point $q\in Q$ and $q$ labeled by $\vee$ or $\wedge$. The rule when the label is $\vee$ is that $q$ gets a black pebble iff there is at least one black pebble on $p_1$ and $p_2$. If the label is $\wedge$ then $q$ gets a white pebble iff there is at least one white pebble on $p_1$ and $p_2$. (So we interpret black pebbles as 1s and white as 0s.)

The nodes of the point-line game will be pairs $(i,a)$, where $i$ is a position on the tape and $a\in A$ is a symbol. They will be connected by an arrow when a transition $(i,a)\to(i+1,b)$ is possible. Each node will be labeled by a player, where we rename Player~I to Black and Player~II to White. The points of a node $(i,a)$ will be all triples $(i,a,c)$, $c\in A$. The lines between points are defined as follows. If $(i,a,c_1)$, $(i,a,c_2)$ and $(i+1,b,d)$ are positions such that Player~I is to decide to move from $(i+1,b,d)$ to either $(i,a,c_1)$ or $(i,a,c_2)$, then $(i+1,b,d)$ is labeled~$\vee$. If for this triple, it is Player~II who is to move, then it is labeled by~$\wedge$. The initial node is $(1,\Lambda)$ and points are the end positions $(1,\Lambda,c)$, $c\in A$.%
\footnote{We assume that the plays always end at the left-most position.} 
The initial position of pebbles consists of black pebbles placed on the winning positions of Player~I, white pebbles are on winning positions of Player~II. 

It is clear that this point-line game simulates the depth 2 game in the sense of general winning strategies. What we must show is that a \emph{positional} strategy in a depth~2 game can be translated into a \emph{positional} strategy in the point-line game. But this is also easy. A positional strategy in a depth~2 game determines, in particular, what a player should do in a position in the first round of the game. These positions correspond to nodes in the point-line game and one can use the same actions.

It remains to show that we can simulate disjunctions and conjunctions by the standard rule of the point-line game. First we observe that we can assume, w.l.o.g., that for any pair of nodes $P$ and $Q$ connected by an arrow $P\to Q$ there is at most one conjunction or disjunction and all other lines are as in the standard game, i.e., they only copy pebbles. This can be achieved by inserting $k$ new nodes if there are $k+1$ conjunctions or disjunctions. 

Suppose we have $(i,a)\to (i+1,b)$ and there are two points $(i,a,c_1),(i,a,c_2)$ in node $(i,a)$ with lines going to one point $(i+1,b,c)$ in node $(i+1,b)$. We now suppose that there can only be one such point $(i+1,b,c)$. Suppose $(i+1,b,c)$ is labeled~$\vee$. Then we insert three new nodes $D$, $D_1$, and $D_2$ and replace the arrow $(i,a)\to (i+1,b)$ with $(i,a)\to D$, $D\to D_1\to (i+1,b)$ and $D\to D_2\to (i+1,b)$. $D$ is labeled Black; the nodes $D_i$ are unlabeled, because there is only one arrow going out of each of them. $D$, $D_1$, and $D_2$ have the same points as $(i,a)$, except that $D_1$ misses $c_2$ and $D_2$ misses $c_1$. The lines between $(i,a)$ and $D$ and between $D$ and $D_1$ and $D_2$, and between $D_1$ and $D_2$ and $(i+1,b)$ connect the corresponding points  except that the lines between the missing points are missing. As a result, Black is able to get a black pebble on $(i+1,b,d)$ iff there is at least one black pebble on $(i,a,c_1)$, or $(i,a,c_2)$. To simulate $\wedge$ we only need to label $D$ by White.

Again, we have to show that the reduction reduces \emph{positional} winning strategies to \emph{positional} winning strategies. Let a strategy for Player~I be given. The translation to the point-line game is straightforward except for the case of the new nodes introduced because of a point labeled~$\vee$. Consider the situation in the previous paragraph. Then we need to define the positional strategy for Black when he is playing at node $D$. For this, we will use his positional strategy when he is playing at position $(i+1,b,c)$ in the depth~2 game: if this strategy is to go to $(i,a,c_j)$, then in the point-line game Black's strategy will be to go from $D$ to $D_j$.

Thus we have shown:

\bpr
Depth 2 game schemas and point-line game schemas strongly polynomially reduce to each other.
\epr

\begin{corollary}
Point-line schemas interpolate $\Pi^s_3$ proofs (equivalently, depth 1 sequent calculus proofs) in the sense of Theorem\ref{t7.1}.
\end{corollary}

\paragraph{General winning strategies.}
We will now present another possible way of viewing point-line game schemas (hence also depth 2 game schemas), but now we will not restrict ourselves to positional winning strategies. Since in every finite game one player has a winning strategy, the point-line game schema defines a \emph{total} monotone Boolean function if we consider all winning strategies.

To motivate what follows, let us first recall how one can view monotone Boolean circuits (hence also depth 1 game schemas). Instead of the standard way where a monotone Boolean circuit is presented as a device computing with bits, one can view it as a way of defining monotone Boolean functions. At a leaf labeled by $x_i$ of the underlying graph we compute the function $x_i$, at a node labeled by $\vee$ we compute the disjunction of the functions defined on the predecessors of the node and similarly on a node labeled $\wedge$. Then the circuit defines the function computed at the root.\footnote{General Boolean circuits can, certainly, be treated in the same way, but in this article we focus on monotone functions.}

Now suppose we are given a point-line game schema. We define functions computed on the nodes of the underlying directed acyclic graph in a similar way, but we will have different variables for every node. We introduce a variable for every point in a node and the function computed at the node will be a function of these variables. If $L$ is a leaf with the unique point $l$, then $f_L$ is the function of one variable $l$ that is the value of this variable. 
Let $P\to Q$, $P\to S$ be arrows in the graph, let $P$ belong to Black, and let $f_Q$ and $f_S$ be functions computed at nodes $Q$ and $S$. Let $p_1\dts p_k\in P$, $q_1\dts q_m\in Q$, and $s_1\dts s_n\in S$ be the points of the three nodes, which we will view as variables of the functions $f_P$, $f_Q$ and $f_S$ respectively. Furthermore, we will view the lines between $P\to Q$ and $P\to S$ as substitutions $\sigma_{PQ}$ and $\sigma_{PS}$, where $\sigma_{PQ}(q_i)=p_j$ if there is a line from $p_j$ to $q_i$, and $\sigma_{PQ}(q_i)=0$ if there is no line from $P$ to $q_i$, and similarly for $\sigma_{PS}$. Then we define
$$f_P(\bar p):=f_Q(\sigma_{PQ}(\bar q))\vee f_R(\sigma_{PS}(\bar s)).$$
If $P$ belongs to White, then the definition is dual ($\vee$ replaced by $\wedge$ and $1$ replaced by $0$). It is not difficult to see that the function at the root $R$ computes who has a winning strategy. We state it as a proposition for further reference.

\bpr
If we interpret black pebbles as 1s and white pebbles as 0s, then the function computed at the root is 1 if Black has a winning strategy and 0 if White has a winning strategy.
\epr

Exponential lower bounds on the size of monotone Boolean circuits of explicitly defined monotone Boolean functions have been proved by the approximation method invented by Razborov~\cite{razb-appr}. In the contemporary presentation this method uses $k$-DNFs and $k$-CNFs, for a suitable $k$ to approximate functions computed at the nodes of the circuit. One shows that (1) at each node only a very small error is introduced and (2) the given function cannot be approximated with a small error. The essence of the method is that the error set of approximating the function computed by the circuit is the union of errors introduced at the nodes, so if the circuit is small this set also has to be small. While the above generalization of monotone Boolean circuits is very similar to the standard monotone Boolean circuits, the approximation method fails in this case. The reason is that the error set introduced at a node is not connected with the set of variables of $f_R$. Various substitutions produce various versions of the error set and thus the total size of the copies can eventually be exponentially larger, even if the circuit has polynomial size.

Note, however, that this computational model may be much stronger than what we need.%
\footnote{We know that every disjoint NP pair can be reduced to the decision who has a winning strategy in a depth 2 game, but it is likely that this decision problem is, in fact, PSPACE complete.} 
We only need computations that tell us who has a \emph{positional} winning strategy; such computation models may be more amenable to lower bounds. 

\subsection{Separation of depth 1 and depth 2 game schemas}

The \emph{clique-coloring tautology} $CC_{m,n,k}$, for $k<m<n$, states that there is no graph on $n$ vertices that has a clique of size $m$ and can be colored by $k$ colors. As an unstatisfiable $CNF$ formula it is formalized by using three sets $M,N,K$, $|M|=m$, $|N|=n$, $|K|=k$, mappings $f:M\to N$ and $g:N\to K$, and a graph $G$ on $N$ and saying the $F$ is one-to-one, $F$ maps $M$ to a clique in $G$ and $g$ is a coloring of $G$. 

The \emph{clique-coloring function} $cc_{m,n,k}$, for $k<m<n$, is the \emph{partial} monotone Boolean function defined on graphs on $n$ vertices that is~1 if the graph has a clique of size $m$, is~0 if the graph is $k$-colorable, and undefined otherwise.

The clique-coloring tautology follows from the \emph{pigeon-hole principle}, because if we compose $f$ with $g$ we get a one-to-one mapping from $M$ to $K$. If $m=2k$, then such a \emph{weak pigeon-hole principle} is provable in depth~1 sequent calculus by proofs of size $n^{(\log n)^{O(1)}}$, which can be used to show that also the clique-coloring tautology $CC_{m,n,k}$ has proofs of asymptotically the same size; see~\cite{krajicek}, Section~18.7. This implies, by our Theorem~\ref{t7.1}, that the clique-coloring function $cc_{m,n,k}$ can be represented by depth~2 game schemas of size  $n^{(\log n)^{O(1)}}$. 

On the other hand, by classical lower bounds on monotone Boolean circuits~\cite{razb-appr,alon-boppana}, any monotone Boolean circuit that computes the clique-coloring function $cc_{m,n,k}$ has exponential size $2^{n^\epsilon}$, $\epsilon>0$, in particular for $m=2k$ and $n=m^3$.

\begin{corollary}
There exists a sequence of partial monotone Boolean functions that can be represented by polynomial size depth~2 schemas (equivalently, point-line schemas) but depth~1 schemas (equivalently, monotone circuits) require exponential size.
\end{corollary}
\bprf
Pad $cc_{m,n,k}$ where $m=2k$ and $n=m^3$ with $n^{(\log n)^{O(1)}}$ dummy bits.
\eprf

%%%%%%%%%%%%%%%%%%%%%%%%%%%%%%%%%%%%%%%%%%%%%%%%%%%%%%%%%%

\section{Open problems}

The next challenge is to characterize the canonical pair of unbounded Frege systems. This pair is polynomially equivalent to the interpolation pair.

\begin{problem}
Characterize the canonical and interpolation pairs of Frege proof systems.
\end{problem}

It seems that our approach should work also in this case. A position in a game obtained from a Symmetric Calculus proof can be determined by a subformula of a formula in the proof. Hence the number of positions in the game is polynomially bounded. The definition of the game should accordingly be modified to allow only polynomial number of positions.

%% But there are other issues. The first one is the following problem, which is also interesting in its own right.

%% \begin{problem}
%% Does the Symmetric Calculus polynomially simulate Frege systems?
%% \end{problem}

A characterization of canonical and interpolation pairs would also be interesting for other weak systems. In particular, the system with disjunctions of \emph{parities} of literals, usually referred to as $Res(Lin)$, is currently intensively studied, but lower bounds have been obtained only for tree-like proofs.%
%\footnote{This system \emph{seems} to be weak, but in fact, we do not know if it can be simulated by a bounded depth Frege system.}

\begin{problem}
Characterize the canonical and interpolation pairs of $Res(Lin)$.
\end{problem}

What we find the most desirable is to extend lower bound methods to stronger computational models. In this article we have presented game schemas as computation models for monotone Boolean functions. The weakest one for which we do not have lower bounds are depth~2 game schemas.

\begin{problem}
Prove a superpolynomial lower bound on depth 2 game schemas representing an explicit monotone partial Boolean function.
\end{problem}

%%%%%%%%%%%%%%%%%%%%%%%%%%%%%%%%%%%%%%%%%%%%%%%%%%%%%%%%%%%%%%%%%%%%%%%

%%%%%%%%%%%%%%%%%%%%%%%%%%%%%%%%%%%%%%%%%%%%%%%%%%%%%%%%%%%%%%%%%%%%%%%

\section*{Appendix}

\subsection*{A1. First order theories and propositional proof systems}

We briefly mention this subject, although we do not use the connection between first order theories and propositional proofs in this article. This section is also an apology why we are not using first order theories.

In fact, we are primarily interested in weak first order theories and study propositional proof complexity because it is a useful tool to prove independence from these theories. The connection first appeared in the seminal article of Stephen Cook~\cite{cook}. A different form was studied by Paris and Wilkie in~\cite{paris-wilkie}. The latter one is more relevant to this work because it connects provability in bounded arithmetic and the length of proofs in bounded depth Frege systems. They extended bounded arithmetic by a new uninterpreted predicate $P$ and added induction for bounded formulas in the extended language. They showed that if a $\Pi^0_1$ sentence is provable in such a theory, then a sequence of tautologies constructed from the sentence has proofs of polynomial lengths in a Frege system restricted to formulas of some constant depth. Using this relation one can show, e.g., that the pigeonhole principle stated with $P$ is not provable in the extended bounded arithmetic.

One can get closer relation between the theories and bounded depth Frege systems if one considers particular fragments. E.g., Buss's theory $T^i_2$ extended to $T^i_2[P]$ leads to \emph{quasipolynomial} Frege proofs in which formulas have depth $i$ with the additional restriction that the \emph{bottom fan-in is polylogarithmic.} In order to get a tight connection, Beckmann et al.~\cite{BPT} introduced special first order theory that capture precisely provability in the depth~$d$ sequent calculus for each $d\geq 0$. Thus one can prove polynomial upper bounds on the lengths of proofs of sequences of tautologies by arguing in a first order theory, which is often more convenient. 

We could have used this connection to prove one part of our result, viz., Lemma~\ref{l5.2}, but we opted not to. In order to use a first order theory, we would have to describe the translation of first order formulas to propositional formulas and eventually the proof would not be much different. The difference would be essentially only in using quantifiers instead of big conjunctions and disjunctions. Furthermore, the theories of~\cite{BPT} are not so well established as $T^i_2[P]$  and we might need to sort out many of details.

Another reason for not using first order theories is to have this article selfcontained.

\subsection*{A2. A remark on simulating cuts}

The simulation of cuts in Lemma~\ref{l-1} is a recursive procedure. It is important that it is run in a ``depth-first'' way. This means that after we split $A\vee(C_1\wedge C_2)$ into  $(A\vee C_1)\wedge(A\vee C_2)$, we first simulate cut with $C_1$ completely and only then we simulate cut with $C_2$. Here is an example.

Consider $(A\vee C)\wedge(B\vee \neg C)$ where $C=C_1\wedge C_2$, $C_1=p\vee q$ and $C_2=r\vee s$. Then the proof will be:
\[

\begin{array}{ccl}
&(A\vee C)\wedge(B\vee \neg C)&\\
\cline{1-2}
&(A\vee C_1)\wedge(A\vee C_2)\wedge(B\vee \neg C_1\vee \neg C_2)&\mbox{\quad by distributivity}\\
\cline{1-2}
&(A\vee p\vee q)\wedge(A\vee C_2)
\wedge(B\vee \neg p\vee \neg C_2)\wedge(B\vee \neg q\vee \neg C_2)&\mbox{\quad by distributivity}
\\
\cline{1-2}
&(A\vee C_2)
\wedge(A\vee B\vee \neg C_2)&\mbox{\quad by resolution with $p$ and $q$}
\\
\cline{1-2}
&(A\vee r\vee s)
\wedge(A\vee B\vee \neg r)
\wedge(A\vee B\vee \neg s)&\mbox{\quad by distributivity}
\\
\cline{1-2}
&A\wedge B&\mbox{\quad by resolution with $r$ and $s$}
\end{array}\]
If we distributed $C_2$ immediately after distributing $C_1$, we would get
\[
\dots\wedge(B\vee \neg p\vee \neg r)
\wedge(B\vee \neg p\vee \neg s)
\wedge(B\vee \neg q\vee \neg r)
\wedge(B\vee \neg q\vee \neg s)
\]
This is like rewriting the DNF $(\neg p\wedge\neg q)\vee(\neg r\wedge \neg s)$ into the CNF 
$
(\neg p\vee \neg r)
\wedge(\neg p\vee \neg s)
\wedge(\neg q\vee \neg r)
\wedge(\neg q\vee \neg s).
$
In general, this operation leads to an exponential blowup.

\bye